\documentclass{gtart}

\def\ifplaintex{\expandafter\ifx\csname documentclass\endcsname\relax}

\ifplaintex 
\else
\fi

\expandafter\ifx\csname beginpicture\endcsname\relax
\expandafter\ifx\csname documentclass\endcsname\relax
\input pictex \else
\input prepictex \input pictex \input postpictex \fi\fi

\def\gt{{\mathsurround=0pt\it $\cal G\mskip-2mu$eometry \&\ 
$\cal T\!\!$opology}}        

\def\gtp{{\mathsurround=0pt\it $\cal G\mskip-2mu$eometry \&\ 
$\cal T\!\!$opology $\cal P\!$ublications}}  

\def\lognumber#1{\def\thelognumber{#1}}
\def\volumenumber#1{\def\thevolumenumber{#1}}
\def\papernumber#1{\def\thepapernumber{#1}}
\def\volumeyear#1{\def\thevolumeyear{#1}}

\def\pagenumbers#1#2{\def\startpage{#1}\def\finishpage{#2}}
\def\published#1{\def\publishdate{#1}}
\def\proposed#1{\def\theproposer{#1}}
\def\seconded#1{\def\theseconders{#1}}
\def\received#1{\def\receiveddate{#1}}
\def\revised#1{\def\reviseddate{#1}}
\def\accepted#1{\def\accepteddate{#1}}
\def\asciititle#1{\def\theasciititle{#1}}

\long\def\asciiabstract#1{\long\def\theasciiabstract{#1}}
\def\asciikeywords#1{\def\theasciikeywords{#1}}

\let\\\par\let\thelognumber\relax
\let\thevolumenumber\relax\let\thepapernumber\relax
\let\thevolumeyear\relax\let\thesamplenumber\relax\let\startpage\relax
\let\finishpage\relax\let\publishdate\relax\let\receiveddate\relax
\let\reviseddate\relax\let\accepteddate\relax\let\theasciititle\relax
\let\theasciiauthors\relax
\let\theasciiabstract\relax\let\theasciikeywords\relax
\let\theasciiemail\relax\let\theshortauthors\relax\let\theshorttitle\relax

\long\def\maketitlep{   

\count0=\startpage

\gt\hfill      
\beginpicture
\setcoordinatesystem units <0.33truein, 0.33truein> point at 2.2 0.9
\setplotsymbol ({$\cal G$})
\plotsymbolspacing=9truept
\circulararc 315 degrees from 0 1 center at 0 0
\setplotsymbol ({$\cal T$})
\circulararc 315 degrees from 1 -1 center at 1 0
\endpicture
\break
{\small\ifx\thesamplenumber\relax 
Volume \else Sample
\fi\thevolumenumber\ (\thevolumeyear)
\startpage--\finishpage\nl
Published: \publishdate}
\vglue 0.5truein plus 0.4fil minus 0.1truein

{\parskip=0pt\leftskip 0pt plus 1fil\def\\{\par\smallskip}{\ifplaintex\large
\else\Large\fi\bf\thetitle}\par\medskip}   

\vglue 0pt plus 0.1fil 

{\parskip=0pt\leftskip 0pt plus 1fil\def\\{\par}{\sc\theauthors}
\par\medskip}

\vglue 0pt plus 0.1fil 

{\small\parskip=0pt\let\newline\\
{\leftskip 0pt plus 1fil\def\\{\par}{\sl\theaddress}\par}
\expandafter\ifx\theemail\relax    
\relax\else\vglue 5pt plus 0.02fil minus 2pt\def\\{\stdspace{\rm 
and}\stdspace} 
\cl{Email:\stdspace\tt\theemail}\fi
\ifx\theurl\relax                  
\relax\else\vglue 5pt plus 0.02fil minus 2pt\def\\{\stdspace{\rm 
and}\stdspace}
\cl{URL:\stdspace\tt\theurl}\fi\par}

\vglue 7pt plus 0.3fil minus 3pt

{\bf Abstract}
\vglue 5pt plus 0.1fil minus 2pt

\theabstract

\vglue 7pt plus 0.3fil minus 3pt

{\bf AMS Classification numbers}\quad Primary:\quad \theprimaryclass

Secondary:\quad \thesecondaryclass

\vglue 5pt plus 0.3fil minus 2pt

{\bf Keywords}\quad \thekeywords

\vglue 10pt plus 0.5fil minus 5pt

{\small  Proposed: \theproposer\hfill Received: \receiveddate\nl
Seconded: \theseconders\hfill 
\ifx\reviseddate\relax                         
Accepted: \accepteddate                        
\else
Revised: \reviseddate                          
\fi}
\eject
}       

\let\maketitlepage\maketitlep
\let\maketitle\maketitlepage

\font\phead=cmsl9 scaled 950
\font=cmsl9 scaled 1050
\font\pnum=cmbx10 scaled 913
\font\lnum=cmbx10 
\font\pfoot=cmsl9 scaled 950
\font=cmsl9 scaled 1050
\ifplaintex
\headline{\vbox to 0pt{\vskip -4.5mm\line{\small\phead\ifnum
\count0=\startpage ISSN 1364-0380 (on line)
1465-3060 (printed) \hfill {\pnum\folio}\else\ifodd\count0\def\\{ }%
\ifx\theshorttitle\relax\thetitle\else\theshorttitle\fi\hfill{\pnum\folio}
\else\def\\{ and }{\pnum\folio}\hfill\ifx\theshortauthors\relax\theauthors
\else\theshortauthors\fi\fi\fi}\vss}}
\footline{\vbox to 0pt{\vglue 0mm\line{\small\pfoot\ifnum\count0=\startpage
\copyright\ \gtp\hfill\else
\gt, Volume \thevolumenumber\ (\thevolumeyear)\hfill\fi}\vss
}}
\else
\makeatletter
\def\@oddhead{{\small\lhead\ifnum\count0=\startpage ISSN 1364-0380 (on line)
1465-3060 (printed) \hfill {\lnum\number\count0}\else\ifodd\count0
\def\\{ }\ifx\theshorttitle\relax \thetitle \else\theshorttitle\fi\hfill
{\lnum\number\count0}\else\def\\{ and }{\lnum\number\count0}
\hfill\ifx\theshortauthors\relax 
\theauthors\else\theshortauthors\fi\fi\fi}}\def\@evenhead{\@oddhead}
\def\@oddfoot{\small\lfoot\ifnum\count0=\startpage\copyright\ \gtp\hfill\else
\gt, Volume \thevolumenumber\ (\thevolumeyear)\hfill\fi}
\def\@evenfoot{\@oddfoot}
\makeatother
\fi

\newwrite\gtoutfile
\long\gdef\makeheadfile{  
{\def\\{, }\def\s{ }
\immediate\openout\gtoutfile head.xxx
\immediate\write\gtoutfile{To: math@arxiv.org}
\immediate\write\gtoutfile{Subject: put or rep NNNNN:pppp}
\immediate\write\gtoutfile{--text follows this line--}
\immediate\write\gtoutfile{Proxy-for: \ifx\theasciiauthors\relax
\theauthors\else\theasciiauthors\fi\s<\ifx\theasciiemail\relax\theemail\else\theasciiemail\fi>}
\immediate\write\gtoutfile{\noexpand\\}
\immediate\write\gtoutfile{Authors: \ifx\theasciiauthors\relax
\theauthors\else\theasciiauthors\fi}
{\def\\{ }\immediate\write\gtoutfile{Title: \ifx\theasciititle\relax
\thetitle\else\theasciititle\fi}}
\immediate\write\gtoutfile{Subj-class: GT or SG or MG etc}
\immediate\write\gtoutfile{MSC-class: \theprimaryclass\ifx\thesecondaryclass\relax\else, \thesecondaryclass\fi}
\immediate\write\gtoutfile{Journal-ref: Geom. Topol. \thevolumenumber
(\thevolumeyear) \startpage-\finishpage}
\immediate\write\gtoutfile{Comments: Published by Geometry and Topology at}
\immediate\write\gtoutfile{\s\s http://www.maths.warwick.ac.uk/gt/GTVol\thevolumenumber/paper\thepapernumber.abs.html}
\immediate\write\gtoutfile{\noexpand\\}
\immediate\write\gtoutfile{}
\ifx\theasciiabstract\relax
\immediate\write\gtoutfile{\theabstract}\else
\immediate\write\gtoutfile{\theasciiabstract}\fi
\immediate\write\gtoutfile{}
\immediate\write\gtoutfile{\noexpand\\}
\immediate\write\gtoutfile{}
\immediate\closeout\gtoutfile}}  

\def\maketitlepage{\maketitlep\makeheadfile}
\let\maketitle\maketitlepage

\def\ifplaintex{\expandafter\ifx\csname documentclass\endcsname\relax}

\ifplaintex 
\else
\fi

\expandafter\ifx\csname beginpicture\endcsname\relax
\expandafter\ifx\csname documentclass\endcsname\relax
\input pictex \else
\input prepictex \input pictex \input postpictex \fi\fi

\def\gt{{\mathsurround=0pt\it $\cal G\mskip-2mu$eometry \&\ 
$\cal T\!\!$opology}}        

\def\gtp{{\mathsurround=0pt\it $\cal G\mskip-2mu$eometry \&\ 
$\cal T\!\!$opology $\cal P\!$ublications}}  

\def\lognumber#1{\def\thelognumber{#1}}
\def\volumenumber#1{\def\thevolumenumber{#1}}
\def\papernumber#1{\def\thepapernumber{#1}}
\def\volumeyear#1{\def\thevolumeyear{#1}}

\def\pagenumbers#1#2{\def\startpage{#1}\def\finishpage{#2}}
\def\published#1{\def\publishdate{#1}}
\def\proposed#1{\def\theproposer{#1}}
\def\seconded#1{\def\theseconders{#1}}
\def\received#1{\def\receiveddate{#1}}
\def\revised#1{\def\reviseddate{#1}}
\def\accepted#1{\def\accepteddate{#1}}
\def\asciititle#1{\def\theasciititle{#1}}

\long\def\asciiabstract#1{\long\def\theasciiabstract{#1}}
\def\asciikeywords#1{\def\theasciikeywords{#1}}

\let\\\par\let\thelognumber\relax
\let\thevolumenumber\relax\let\thepapernumber\relax
\let\thevolumeyear\relax\let\thesamplenumber\relax\let\startpage\relax
\let\finishpage\relax\let\publishdate\relax\let\receiveddate\relax
\let\reviseddate\relax\let\accepteddate\relax\let\theasciititle\relax
\let\theasciiauthors\relax
\let\theasciiabstract\relax\let\theasciikeywords\relax
\let\theasciiemail\relax\let\theshortauthors\relax\let\theshorttitle\relax

\long\def\maketitlep{   

\count0=\startpage

\gt\hfill      
\beginpicture
\setcoordinatesystem units <0.33truein, 0.33truein> point at 2.2 0.9
\setplotsymbol ({$\cal G$})
\plotsymbolspacing=9truept
\circulararc 315 degrees from 0 1 center at 0 0
\setplotsymbol ({$\cal T$})
\circulararc 315 degrees from 1 -1 center at 1 0
\endpicture
\break
{\small\ifx\thesamplenumber\relax 
Volume \else Sample
\fi\thevolumenumber\ (\thevolumeyear)
\startpage--\finishpage\nl
Published: \publishdate}
\vglue 0.5truein plus 0.4fil minus 0.1truein

{\parskip=0pt\leftskip 0pt plus 1fil\def\\{\par\smallskip}{\ifplaintex\large
\else\Large\fi\bf\thetitle}\par\medskip}   

\vglue 0pt plus 0.1fil 

{\parskip=0pt\leftskip 0pt plus 1fil\def\\{\par}{\sc\theauthors}
\par\medskip}

\vglue 0pt plus 0.1fil 

{\small\parskip=0pt\let\newline\\
{\leftskip 0pt plus 1fil\def\\{\par}{\sl\theaddress}\par}
\expandafter\ifx\theemail\relax    
\relax\else\vglue 5pt plus 0.02fil minus 2pt\def\\{\stdspace{\rm 
and}\stdspace} 
\cl{Email:\stdspace\tt\theemail}\fi
\ifx\theurl\relax                  
\relax\else\vglue 5pt plus 0.02fil minus 2pt\def\\{\stdspace{\rm 
and}\stdspace}
\cl{URL:\stdspace\tt\theurl}\fi\par}

\vglue 7pt plus 0.3fil minus 3pt

{\bf Abstract}
\vglue 5pt plus 0.1fil minus 2pt

\theabstract

\vglue 7pt plus 0.3fil minus 3pt

{\bf AMS Classification numbers}\quad Primary:\quad \theprimaryclass

Secondary:\quad \thesecondaryclass

\vglue 5pt plus 0.3fil minus 2pt

{\bf Keywords}\quad \thekeywords

\vglue 10pt plus 0.5fil minus 5pt

{\small  Proposed: \theproposer\hfill Received: \receiveddate\nl
Seconded: \theseconders\hfill 
\ifx\reviseddate\relax                         
Accepted: \accepteddate                        
\else
Revised: \reviseddate                          
\fi}
\eject
}       

\let\maketitlepage\maketitlep
\let\maketitle\maketitlepage

\font\phead=cmsl9 scaled 950
\font=cmsl9 scaled 1050
\font\pnum=cmbx10 scaled 913
\font\lnum=cmbx10 
\font\pfoot=cmsl9 scaled 950
\font=cmsl9 scaled 1050
\ifplaintex
\headline{\vbox to 0pt{\vskip -4.5mm\line{\small\phead\ifnum
\count0=\startpage ISSN 1364-0380 (on line)
1465-3060 (printed) \hfill {\pnum\folio}\else\ifodd\count0\def\\{ }%
\ifx\theshorttitle\relax\thetitle\else\theshorttitle\fi\hfill{\pnum\folio}
\else\def\\{ and }{\pnum\folio}\hfill\ifx\theshortauthors\relax\theauthors
\else\theshortauthors\fi\fi\fi}\vss}}
\footline{\vbox to 0pt{\vglue 0mm\line{\small\pfoot\ifnum\count0=\startpage
\copyright\ \gtp\hfill\else
\gt, Volume \thevolumenumber\ (\thevolumeyear)\hfill\fi}\vss
}}
\else
\makeatletter
\def\@oddhead{{\small\lhead\ifnum\count0=\startpage ISSN 1364-0380 (on line)
1465-3060 (printed) \hfill {\lnum\number\count0}\else\ifodd\count0
\def\\{ }\ifx\theshorttitle\relax \thetitle \else\theshorttitle\fi\hfill
{\lnum\number\count0}\else\def\\{ and }{\lnum\number\count0}
\hfill\ifx\theshortauthors\relax 
\theauthors\else\theshortauthors\fi\fi\fi}}\def\@evenhead{\@oddhead}
\def\@oddfoot{\small\lfoot\ifnum\count0=\startpage\copyright\ \gtp\hfill\else
\gt, Volume \thevolumenumber\ (\thevolumeyear)\hfill\fi}
\def\@evenfoot{\@oddfoot}
\makeatother
\fi

\newwrite\gtoutfile
\long\gdef\makeheadfile{  
{\def\\{, }\def\s{ }
\immediate\openout\gtoutfile head.xxx
\immediate\write\gtoutfile{To: math@arxiv.org}
\immediate\write\gtoutfile{Subject: put or rep NNNNN:pppp}
\immediate\write\gtoutfile{--text follows this line--}
\immediate\write\gtoutfile{Proxy-for: \ifx\theasciiauthors\relax
\theauthors\else\theasciiauthors\fi\s<\ifx\theasciiemail\relax\theemail\else\theasciiemail\fi>}
\immediate\write\gtoutfile{\noexpand\\}
\immediate\write\gtoutfile{Authors: \ifx\theasciiauthors\relax
\theauthors\else\theasciiauthors\fi}
{\def\\{ }\immediate\write\gtoutfile{Title: \ifx\theasciititle\relax
\thetitle\else\theasciititle\fi}}
\immediate\write\gtoutfile{Subj-class: GT or SG or MG etc}
\immediate\write\gtoutfile{MSC-class: \theprimaryclass\ifx\thesecondaryclass\relax\else, \thesecondaryclass\fi}
\immediate\write\gtoutfile{Journal-ref: Geom. Topol. \thevolumenumber
(\thevolumeyear) \startpage-\finishpage}
\immediate\write\gtoutfile{Comments: Published by Geometry and Topology at}
\immediate\write\gtoutfile{\s\s http://www.maths.warwick.ac.uk/gt/GTVol\thevolumenumber/paper\thepapernumber.abs.html}
\immediate\write\gtoutfile{\noexpand\\}
\immediate\write\gtoutfile{}
\ifx\theasciiabstract\relax
\immediate\write\gtoutfile{\theabstract}\else
\immediate\write\gtoutfile{\theasciiabstract}\fi
\immediate\write\gtoutfile{}
\immediate\write\gtoutfile{\noexpand\\}
\immediate\write\gtoutfile{}
\immediate\closeout\gtoutfile}}  

\def\maketitlepage{\maketitlep\makeheadfile}
\let\maketitle\maketitlepage

\lognumber{206}
\volumenumber{5}\papernumber{28}\volumeyear{2001}
\pagenumbers{895}{924}
\received{1 September 2001}
\accepted{31 December 2001}
\proposed{Ronald Fintushel}
\seconded{Robion Kirby, John Morgan}
\revised{2 January 2002}
\published{3 January 2002}

\usepackage{amsmath,amssymb,graphicx}

\def\S{section }

\theoremstyle{plain} \newtheorem{theorem}{Theorem}[section]
\newtheorem{lemma}[theorem]{Lemma}
\newtheorem{corollary}[theorem]{Corollary}
\newtheorem{proposition}[theorem]{Proposition}

\theoremstyle{remark} 
\newtheorem*{notation}{Notation} 
\newtheorem*{claim}{Claim}
\newtheorem*{ack}{Acknowledgement}

\theoremstyle{definition} \newtheorem{defn}[theorem]{Definition}

\def \R {\mathbf{R}}
\def \Z {\mathbf{Z}}

\def\I{\mathbf{I}}
\def\cp{\mathbf{CP}}
\def\cpbar{\overline{\mathbf{CP}}^{2}}
\def\HH{\mathcal{H}}
\def\WW{\mathcal{W}}

\def\OO{\mathcal{O}}

\def\mtilde{\tilde{\mathcal{M}}}
\def\mttilde{\Tilde{\Tilde{\mathcal{M}}}}
\def\ttimes{\tilde{\times}}

 \def\MM{\mathcal{M}}
 \def\mpsc{\mathcal{M}^{+}}
 
 \def\diff{\mathit{Diff}}
\def\SC{\ifmmode{\text{SPIN}^c}\else{$\text{SPIN}^c$}\fi}
\DeclareMathOperator{\SO}{SO}
\DeclareMathOperator{\met}{Met}
\DeclareMathOperator{\psc}{PSC}
\DeclareMathOperator{\spinc}{Spin^c}
\DeclareMathOperator{\pinc}{Pin^c}
\DeclareMathOperator{\pinp}{Pin^+}

\DeclareMathOperator{\pinpm}{Pin^\pm}
\DeclareMathOperator{\coker}{coker}
\def\sw{Seiberg--Witten}
\def\SW{\ifmmode{\text{SW}}\else{$\text{SW}$}\fi}
\def\SWW{\ifmmode{{\text{SW}}_{tot}}\else{${\text{SW}}_{tot}$}\fi}
\def\mfd{$4$--manifold}

\begin{document}

\title[Positive scalar curvature]{Positive scalar curvature,
diffeomorphisms\\and the Seiberg--Witten invariants}
\asciititle{Positive scalar curvature, diffeomorphisms and the
Seiberg-Witten invariants}

\author{Daniel Ruberman}
\address{Department of Mathematics, MS 050\\Brandeis
University\\Waltham, MA 02454-9110, USA}
\email{ruberman@brandeis.edu}

\begin{abstract}
We study the space of positive scalar curvature (psc) metrics on a
4--manifold, and give
examples of simply connected manifolds for which it is disconnected. These
examples imply that concordance of psc metrics does not imply isotopy of such
metrics. This is demonstrated using a modification of the 1--parameter
Seiberg--Witten invariants which we introduced in earlier work. The invariant
shows that the diffeomorphism group of the underlying 4--manifold is
disconnected. We also study the moduli space of positive scalar curvature
metrics modulo diffeomorphism, and give examples to show that this space can be
disconnected. The (non-orientable) 4--manifolds in this case are explicitly
described, and the components in the moduli space are distinguished by a
$Pin^c$ eta invariant.
\end{abstract}

\asciiabstract{
We study the space of positive scalar curvature (psc) metrics on a
4-manifold, and give
examples of simply connected manifolds for which it is disconnected. These
examples imply that concordance of psc metrics does not imply isotopy of such
metrics. This is demonstrated using a modification of the 1-parameter
Seiberg-Witten invariants which we introduced in earlier work. The invariant
shows that the diffeomorphism group of the underlying 4-manifold is
disconnected. We also study the moduli space of positive scalar curvature
metrics modulo diffeomorphism, and give examples to show that this space can be
disconnected. The (non-orientable) 4--manifolds in this case are explicitly
described, and the components in the moduli space are distinguished by a
Pin^c eta invariant.}

\asciikeywords{Positive scalar curvature, Seiberg-Witten equations, isotopy}
\maketitlepage

\section{Introduction}\label{sec:intro}
One of the striking initial applications of the Seiberg--Witten
invariants was to give new obstructions to the existence of Riemannian
metrics of positive scalar curvature on $4$--manifolds.  The vanishing
of the Seiberg--Witten invariants of a manifold admitting such a
metric may be viewed as a non-linear generalization of the classic
conditions~\cite{lichnerowicz:spinors,lawson-michelson,rosenberg-stolz:psc}
derived from the Dirac operator.  If a manifold $Y$ has a metric of
positive scalar curvature, it is natural to investigate the topology
of the space $\psc(Y)$ of all such metrics.  Perhaps the simplest
question which one can ask is whether $\psc(Y)$ is connected; examples
of manifolds
for which it is disconnected were previously known in all dimensions
greater than $4$.  This phenomenon is detected via the index theory of
the Dirac operator, often in conjunction with the
Atiyah--Patodi--Singer index theorem~\cite{aps:I}.

In the first part of this paper, we use a variation of the
$1$--parameter \sw\ invariant introduced in ~\cite{ruberman:isotopy}
to prove that on a simply-connected $4$--manifold $Y$, $\psc(Y)$ can
be disconnected.  Our examples cannot be detected by index theory
alone, ie without the intervention of the Seiberg--Witten equations.
These examples also yield a negative answer, in dimension $4$, to the
question of whether metrics of positive scalar curvature which are
concordant are necessarily isotopic.  Apparently
(cf the discussion in~\cite[\S3 and \S6]{rosenberg-stolz:psc}) this is the
first result of this sort in any dimension other than $2$.

An {\sl a priori} more difficult problem than showing that $\psc$ is
not necessarily connected is to find manifolds for
which the ``moduli space'' $\psc/\diff$ is disconnected.  (The action
of the diffeomorphism group on the space of metrics is by pull-back,
and preserves the subset of positive scalar curvature metrics.)  The
metrics lying in different components of $\psc(Y)$ constructed in the
first part of the paper are obtained by pulling back a positive scalar
curvature metric via one of the diffeomorphisms introduced
in~\cite{ruberman:isotopy}, and hence give no information about
$\psc/\diff$.  Building on constructions of Gilkey~\cite{gilkey:even}
we give explicit examples of non-orientable \mfd s for which the
moduli space is disconnected.  These examples are detected, as
in~\cite{gilkey:even,botvinnik-gilkey:spaceforms}, by an
$\eta$--invariant associated to a $\pinc$ Dirac operator.

The \sw\ invariant for diffeomorphisms introduced in
~\cite{ruberman:isotopy} is not very amenable to calculation, for
reasons explained below.  Hence, in the first section of this paper,
we will give a modification of that construction which yields a more
computable invariant.  From this modified \sw\ invariant we will
deduce the non-triviality of the isotopy class of the diffeomorphisms
described in ~\cite{ruberman:isotopy}, the non-connectedness
of $\psc$, and the fact that concordance does not imply isotopy.  The
modification has the disadvantage of not behaving very
sensibly under composition of diffeomorphisms, but it suffices for the
application to the topology of $\psc$.
\begin{ack}
I would like to thank Alex Suciu for reminding me of Pao's work on
spun lens spaces, as well as Selman Akbulut and Ian Hambleton for
telling me about their approaches to these manifolds. I also thank
Peter Gilkey for a helpful communication and the referee for some 
perceptive comments on the exposition.  The author was partially 
supported
by NSF Grant 9971802.  This paper was completed during a visit to the IHES and
the University of Paris, Orsay, under the partial support of the CNRS.
\end{ack}

\section{Seiberg--Witten invariants of diffeomorphisms}\label{sec:sw}
Let us briefly describe the approach taken in~\cite{ruberman:isotopy}
to defining gauge-theoretic invariants of diffeomorphisms, and
explain a modification of that approach which renders the \sw\ version
of those invariants more usable in the current paper.
\begin{notation}
            A $\spinc$ structure on a manifold $Z$ will be indicated by
            $\Gamma$, representing the $\spinc$ bundles $W^\pm$ and the
            Clifford multiplication $\Lambda^{*}(Z) \times W^\pm \to
            W^\mp$.  The set of perturbations for the \sw\ equations will
            be denoted by:
            $$
            \Pi = \{(g,\mu)\in \met(Z) \times \Omega^2(Z) \, | \, *_g \mu
            =\mu\}
            $$
            We will sometimes identify $\met(Z)$ with its image $\{(g,0)|\,
            g \in \met(Z)\}$, and so will write $g \in \Pi$ instead of
            $(g,0) \in \Pi$.  A Riemannian metric $g$ and a $U(1)$ 
connection $A$
on $L =
\det(W^+)$ determine the Dirac operator $D_A\co \Gamma(W^+) \to
\Gamma(W^-)$.  For
$(g,\mu)
\in \Pi$, the perturbed \sw\ equations for $A$ and $\varphi\in
\Gamma(W^+)$ are written:
\begin{equation}
\begin{cases} D_A\varphi = 0\\ F_A^+ + \imath \mu = q(\varphi)
\end{cases}
\end{equation} The moduli space of solutions to the perturbed \sw\
equations (modulo gauge
equivalence) associated to $h\in \Pi$ will be denoted $\MM(\Gamma;h)$.
\end{notation}
In principle, a $\spinc$ structure involves a choice of Riemannian
metric on $Z$, but as explained in~\cite[\S 2.2]{ruberman:isotopy} a
$\spinc$ structure using one metric gives rise canonically to a
$\spinc$ structure for any metric.

Choose a homology orientation for $Z$, which is the data needed to
orient the determinant line bundle associated to the \sw\ equations
for every choice of $\Gamma$.  Suppose that $\Gamma$ is a $\spinc$
structure on $Z$, such that the moduli space has formal dimension
equal to $-1$.  Thus, for a generic $h_0 \in \Pi$, the moduli space
$\MM(\Gamma;h_{0})$ will be empty.  For a path $h_t \subset \Pi$,
define the $1$--parameter moduli space to be
$$
\mtilde(\Gamma; h) = \bigcup_{t\in [0,1]} \MM(\Gamma;h_t).
$$
If $h$ is generic, then this space is $0$--dimensional, and we can
count its points with signs.  It is shown in~\cite{ruberman:isotopy}
that $\SW(\Gamma;h)$ only depends on the endpoints $h_0$ and $h_1$,
and not on the choice of path.

\begin{defn}\label{1param}
            For generic $h_0,h_1 \in \Pi$, and for a generic path $h\co
            [0,1] \to \Pi$, we will denote by $\SW(\Gamma;h) =
            \SW(\Gamma;h_0, h_1)$ the algebraic count $\#\mtilde(\Gamma;
            h)$.
\end{defn}

The applications of Seiberg--Witten theory to the topology of $\psc$
in this paper depend on the following simple observation.
\begin{lemma}\label{psc}
            Suppose $b_2^+(Y) \ge 2$, and that $g_0, g_1$ are generic
            Riemannian metrics in the same component of $\psc(Y)$. Let
            $h_i = (g_i,\mu_i) \in \Pi$ be generic.  Then
            $\SW(\Gamma;h_0, h_1) = 0$ if the $\mu_i$ are sufficiently
            small.
\end{lemma}

Suppose that $f\co Z \to Z$ is a diffeomorphism satisfying
\begin{equation}\label{preserve}
            f^*(\Gamma) \cong \Gamma\ \mathrm{and} \ \alpha(f) = 1.
\end{equation}
Here $\alpha(f)= \pm 1$ where the sign is that of the determinant of
the map which $f$ induces on $H^2_+(Z)$.  Informally, $\alpha(f) = 1$
means that $f$ preserves the orientation of the moduli space
$\MM(\Gamma)$.  It is shown in~\cite{ruberman:isotopy} that the
quantity $\SW(\Gamma;h_0, f^*(h_0))$ does not depend on the choice of
$h_0$.  In particular, if $Y$ has a metric $g_0$ of positive scalar
curvature, and a diffeomorphism satisfying~\eqref{preserve} has
$\SW(\Gamma;g_0, f^*(g_0)) \neq 0$, then $g_0$ and $ f^*(g_0)$ are in
different components of $\psc$.  It also follows that $f$ is not
isotopic to the identity, although that is not our main concern in
this paper.

In~\cite{ruberman:isotopy}, certain diffeomorphisms were constructed
that have nontrivial Donaldson invariants, but we were unable to
evaluate their Seiberg--Witten invariants.  The reason for this is that
these diffeomorphisms are constructed as compositions of simpler
diffeomorphisms, which do not satisfy~\eqref{preserve}.  We will now
describe how to modify Definition~\ref{1param} so that the resulting
invariants are defined and computable for a broader class of
diffeomorphisms, including those discussed in ~\cite{ruberman:isotopy}.

Let $\SC(Y)$ denote the set of $\spinc$ structures on a manifold $Y$.
For an orientation-preserving diffeomorphism $f\co Y \to Y$, and a
$\spinc$ structure $\Gamma$, let $\OO(f,\Gamma)$ be the orbit of
$\Gamma$ under the natural action of $f$ on \SC\ via pullback.
\begin{defn}\label{spincsum}
            Suppose that $\Gamma$ is a $\spinc$ structure on Y such that
            the Seiberg--Witten moduli space $\MM(\Gamma)$ has formal
            dimension $-1$.  For an arbitrary generic point $h_{0}\in
            \Pi$, define
\begin{equation}
            \SWW(f,\Gamma) = \sum_{\Gamma' \in\
\OO(f,\Gamma)}\SW(\Gamma';h_{0},f^{*}h_{0}). \label{eq:average}
\end{equation}
\end{defn}
If the orbit $\OO(f,\Gamma)$ is finite (for instance if $f$ preserves
the $\spinc$ structure $\Gamma$) then the sum in
Definition~\ref{spincsum} makes evident sense.  More generally, we
have the following consequence of the basic analytical properties of
the \sw\ equations.
\begin{proposition}\label{prop:finite}
           For a generic path $h\co [0,1] \to \Pi$, there are a finite
           number of $\spinc$ structures on $Y$ for which $\MM(Y;h) $ is
           non-empty.
\end{proposition}
\begin{proof}  Since $b_+^2 >2$,  no reducible solutions will arise in
           a generic path of metrics. In a neighborhood of a given path
           $(g,\mu)$ in $\Pi$, there is a uniform lower bound on the
           scalar curvature and the norm of $\mu$.  By a basic argument
           in \sw\ theory, this implies that there are a finite number of
           $\spinc$ structures for which the parameterized moduli space
           has dimension $\ge 0$.  But for generic paths, a moduli space
           with negative dimension will be empty.
\end{proof}

The main point we will need to verify concerning
Definition~\ref{spincsum} is that the quantity defined
in~\eqref{eq:average} is independent of the choice of $h_{0}\in \Pi$.
In discussing the invariance properties of $\SWW(f,\Gamma)$ it is
useful to rewrite it as a sum over a long path in $\Pi$.  For $n \in
\Z$, we write $f_{n}$ for the $n$--fold composition of $f$, with the
understanding that $f_{0}=id$ and $f_{-n} = f_{n}^{-1}$.
Correspondingly, given $h_{0}\in \Pi$, we write $h_{n} =
f_{n}^{*}h_{0}$.  If $h\co [0,1] \to \Pi $ is a path from $h_{0}$ to
$h_{1}$, then there are obvious paths $f_{n-1}^{*}h$ between $h_{n-1}$
and $h_{n}$, which fit together to give a continuous map $\R \to \Pi$.
We distinguish two cases, according to the way in which $f$ acts on
the $\spinc$ structures.
\begin{enumerate}
            \item If $\OO(f,\Gamma)$ is finite, then we can write
            $$ \SWW(f,\Gamma) = \sum_{n=1}^{N}\SW(f_{-n}^{*}\Gamma;h)
            $$
            where $N$ is the smallest positive integer with
            $f_{N}^{*}\Gamma \cong \Gamma$.  \item If $\OO(f,\Gamma)$ is
            infinite, then we can write
            $$ \SWW(f,\Gamma) =
\sum_{n=-\infty}^{\infty}\SW(f_{-n}^{*}\Gamma;h)
            $$
            keeping in mind that the sum has only finitely many non-zero
            terms, for a generic path $h$.
\end{enumerate}
We will usually write $\sum_{n}\SW(f_{-n}^{*}\Gamma;h) $ so as to be
able to to discuss the two cases simultaneously; the reason for using
$f_{-n}^{*}$ instead of $f_{n}^{*}$ should become clear momentarily.
\begin{lemma}\label{lem:longpath}
            If $\alpha(f) = 1$, then $\SWW(f,\Gamma)$ may be rewritten
            $$
            \SWW(f,\Gamma) = \sum_{n}\SW(\Gamma;f_{n}^{*}h) =
            \sum_{n}\SW(\Gamma;f_{n}^{*}h_{0},f_{n+1}^{*}h_{0})
            $$
\end{lemma}
\begin{proof}
            The proof rests on the isomorphism of moduli spaces
            $$
            f^{*}\co \MM(\Gamma;h) \to \MM(f^{*}\Gamma;f^{*}h)
            $$
            given by pulling back the spin bundle, spinors, and
            connection.  This isomorphism preserves orientation precisely
            when $\alpha(f) = 1$.  Note that $\alpha(f_n) = \alpha(f)^n$.
            Thus, if $h$ is either a generic element of $\Pi$ or a
            generic path in $\Pi$, we have a diffeomorphism
            $$
            f_{n}^{*}\co \MM(f_{-n}^{*}\Gamma;h) \to \MM(\Gamma;f_{n}^{*}h)
            $$
            and the lemma follows by summing up over $\OO(f,\Gamma)$.
\end{proof}
\begin{theorem}\label{welldef}
            Suppose that $b_+^2(Y) > 2$ and that $f$ is an orientation
            preserving diffeomorphism with $\alpha(f) = 1$.  Let $h_0,
            k_0 \in \Pi$ be generic, and let $h,k\co [0,1] \to \Pi$ be
            generic paths with $h_1=f^*h_0$ and $k_1 = f^*k_0$.  Then the
            sums defining $\SWW(f,\Gamma)$ using $h_{0}$ and $k_{0}$ are
            equal:
            $$
            \sum_{n}\SW(\Gamma;f_{n}^{*}h) =
            \sum_{n}\SW(\Gamma;f_{n}^{*}k)
            $$
\end{theorem}
\begin{proof}[Proof of Theorem~\ref{welldef}]
            As remarked before, each term in
            $$
            \sum_{n}\SW(\Gamma;f_{n}^{*}h)
            $$
            depends only on the endpoints $h_{n},h_{n+1}$ of the path.
            Start by choosing a generic path $K_{0,t} \subset \Pi$ from
            $h_0 = K_{0,0}$ to $k_0 = K_{0,1}$, and note that $K_{1,t} =
            f^*K_{0,t}$ is necessarily generic for $f^{*}\Gamma$.  Now
            take generic paths $h$ and $k$, and define $K_{s,0} = h_s$
            and $K_{s,1} = k_s$.  Thus we have a well-defined loop in the
            contractible space $\Pi$, which may be filled in with a
            generic $2$--parameter family $K_{s,t}$.

            Consider the $2$--parameter moduli space
            $$
            \mttilde(\Gamma;K) = \bigcup_{s,t} \MM(\Gamma;K_{s,t})
            $$
            which is a smooth compact $1$--manifold with boundary.  (For
            compactness, we need to ensure that there are no reducibles
            in a $2$--parameter family, which is why we require that $b_+^2
            > 2$.)  Any boundary point lies in the interior of one of the
            sides of the rectangle indicated in Figure 1.  Now form the
            union
            $$
            \bigcup_{n} \MM(\Gamma;f_{n}^{*}K)
            $$
            where the union is taken over the same set of $n$ as is the
            sum defining $\SW(f;\Gamma)$.  Note that by construction the
            moduli space corresponding to the left side of the $n^{th}$
            square (ie, the $1$--parameter moduli space corresponding to
            $f_{n}^{*}K_{0,t}$) matches up with the right side of the
            $(n-1)^{st}$ square (corresponding to $f_{n-1}^{*}K_{1,t}$).
            The argument diverges slightly according to whether
            $\OO(f,\Gamma)$ is finite or infinite.  In the former case,
            assume that the orbit has exactly $N$ elements.  Hence we
            have a diffeomorphism from the right hand boundary of the
            $N^{th}$ square \begin{equation} f_{-N}^{*}\co
            \MM(\Gamma;\{f_{N}^{*}K_{0,t}\}) \to
            \MM(f_{-N}^{*}\Gamma;\{K_{0,t}\}) \cong
            \MM(\Gamma;\{K_{0,t}\}) \end{equation} to the left hand
            boundary of the first square.  Note that the last
            isomorphism, because it is defined up to gauge equivalence,
            is independent of the choice of isomorphism between $\Gamma$
            and $f^{*}_{N}\Gamma$.  Moreover, by the assumption that
            $\alpha(f) =1$, this isomorphism is orientation preserving.
            All of the other contributions from the side boundary
            components cancel, according to the observation in the last
            paragraph.  So, counting with signs, the number of points on
            the top boundary of
            $$
            \bigcup_{n} \MM(\Gamma;f_{n}^{*}K)
            $$
            is the same as the number of points on the bottom part of the
            boundary.  It follows that the conclusion of
            theorem~\ref{welldef} holds in the case that $\OO(f,\Gamma)$
            is finite.  When $\OO(f,\Gamma)$ is infinite, we use the
            principle (Proposition~\ref{prop:finite}) that for $|n| >$
            some $N_{0}$, the moduli spaces $\MM(\Gamma;f_{n}^{*}K)$ are
            all empty.  Thus the union of all of the $2$--parameter
            moduli spaces provides a compact cobordism between the union
            of $1$--parameter moduli spaces corresponding to the initial
            point $h_{0} \in \Pi$ and that corresponding to $k_{0}$.
            \end{proof}

\begin{figure}[ht!]\small
\begin{center}
$\mathcal{M}(\Gamma;f_{-1}^*K_{s,t})$\qquad\qquad\qquad$\mathcal{M}
(\Gamma;f_{0}^*K_{s,t})$\qquad\qquad\qquad$\mathcal{M}(\Gamma;f_{+1}^*K_{s,t})$
\\
\smallskip
\includegraphics[scale=.864]{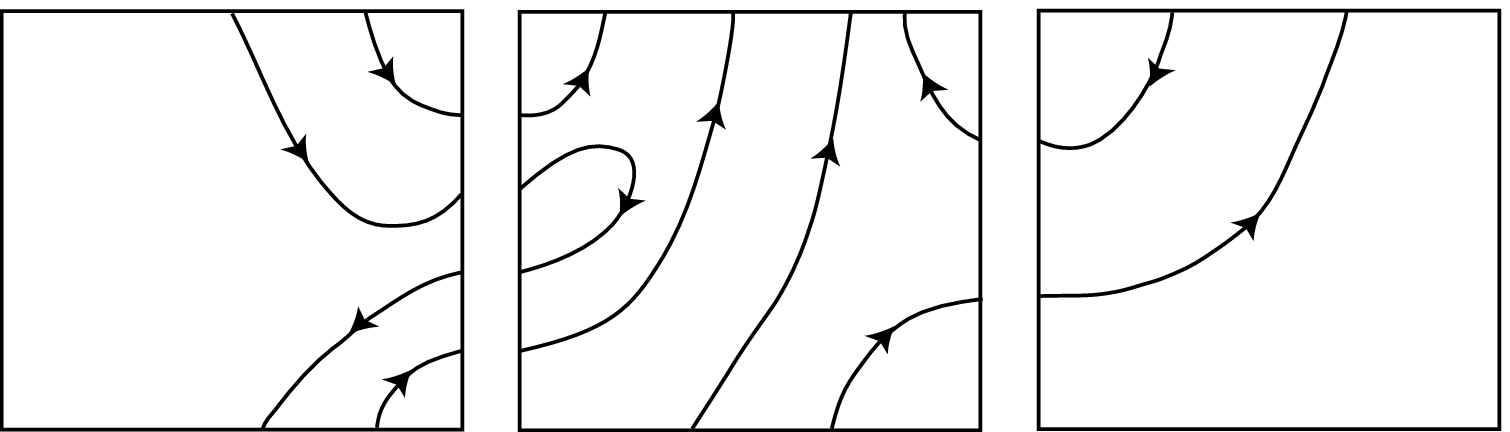} 
\smallskip
Figure 1
\end{center}
\end{figure}

The proof is illustrated above in Figure 1, in the case that
$\OO(f,\Gamma)$ is infinite, and that $\MM(\Gamma;f_{n}^{*}K) \neq
\emptyset$ for $n = -1,0,1$.  If $\OO(f,\Gamma)$ is finite, say with
$N$ elements, there might be components of $\MM(\Gamma;f_{N}^{*}K)$
going off the right-hand edge of the last box, but these would be
matched up with the left-hand edge of the first box.  In either case,
the algebraic number of points on the top and bottom would be the
same.

\section{Basic properties}\label{sec:properties}

\begin{theorem}\label{vanish}
           Suppose $b_2^+(Y) \ge 2$.  If $g$ is a generic Riemannian
           metric of positive scalar curvature, and there is a path in
           $\psc(Y)$ from $g_0$ to $f^*g_0$, then $\SWW(f,\Gamma) = 0$.
\end{theorem}
\begin{proof}
           Since $b_2^+(Y) \ge 2$, there are no reducible solutions to
           the Seiberg--Witten equations for generic metrics or paths of
           metrics. Because $\psc(Y)$ is open in the space of metrics,
           any path may be perturbed to be generic while staying in
           $\psc(Y)$.  Since $b_2^+(Y) \ge 2$, there are no reducible
           solutions to the Seiberg--Witten equations for generic metrics
           or paths of metrics.  Now if $g \subset \psc$ is a path from
           $g_0$ to $f^*g_0$, then all of the paths $f_k^*(g)$ lie in
           $\psc$.  By the standard vanishing theorem, all of the moduli
           spaces
           $$
           \mtilde(\Gamma,f_k^*(g))
           $$
           are empty, and so the invariant $\SWW(f,\Gamma)$ must vanish.
\end{proof}

The isotopy invariance of $\SWW(f,\Gamma)$ is a little more
complicated to prove.  It depends on the device of passing back and
forth between the definition of $\SWW$ as a sum over $\spinc$
structures and its definition as a sum over paths in $\Pi$.

\begin{theorem}\label{isotopy}
           Suppose that $f$ and $g$ are isotopic diffeomorphisms of $Y$
           for which $\SWW(f,\Gamma)$ are defined.  Then $\SWW(f,\Gamma)
           = \SWW(g,\Gamma)$.
\end{theorem}
\begin{proof}
           In the proof we use the notation $fg$ for the composition
           $f\circ g$ of functions, and $h\cdot k$ for the composition of
           paths (the path $h$ followed by the path $k$.)  We will use an
           isotopy between $f$ and $g$ to construct a diffeomorphism
           between the various $1$--parameter moduli spaces, taking
           advantage of the fact that any with the correct endpoints may
           be used to calculate $\SW(\Gamma;g_n^*h_0,g_{n+1}^*h_0)$.
           First choose an isotopy $F_t\co  Y \to Y$ such that $F_0 = id$
           and $F_1 = f_{-1}g$.  Let $h_0 \in \Pi$ be generic, and choose
           a generic path $h$ from $h_0$ to $f^*h_0$.  By definition, $$
           \SWW(f,\Gamma) = \sum_n \SW(\Gamma;h_n,h_{n+1}) $$ where $h_n
           = f_n^*h_0$.  We use the path $f_n^*h$ to calculate the
           $n^{th}$ term in the sum.

           Let
           $$
           F_{k,t} = f_{-k}F_t g_k,\ t \in [0,1]
           $$
           which gives an isotopy of $ f_{-k} g_k$ to $ f_{-(k+1)}
           g_{k+1}$.  We calculate $\SWW(g,\Gamma)$ using $h_0$ for the
           initial point.  The path $F_{k,t}^* f_k^* h_t = g_k^*F_t^*
           h_t$ has endpoints
           $$
           g_k^*F_0^* h_0 = g_k^*h_0
           $$
           for $t=0$ and
           $$
           g_k^*F_1^* h_1 = g_k^* g^* f_{-1}^* f^* h_0 = g_{k+1}^*h_0
           $$
           for $t=1$.  So the $1$--parameter moduli spaces
           $$
           \MM(\Gamma;\{F_{k,t}^* f_k^* h_t\})
           $$
           may be used to calculate $\SWW(g,\Gamma)$.

           But $F_{k,t}$ is homotopic to the identity, and so acts as the
           identity on the set of $\spinc$ structures on $Y$.  Hence
           there is an orientation preserving diffeomorphism
           $$
           F_{k,t}^*\co  \MM(\Gamma;F_{k,t}^* f_k^* h_t) \to
\MM((F_{k,t}^{-1})^* \Gamma; f_k^* h_t) \cong \MM(\Gamma;f_k^* h_t)
           $$
           which is the moduli space used to calculate $\SWW(f;\Gamma)$.
           (Note that the gen-\break ericity of the path $F_{k,t}^* f_k^* h_t$
           follows from this diffeomorphism.)  Therefore\break $\SWW(f;\Gamma)
           = \SWW(g,\Gamma)$.
\end{proof}

One thing which would be useful to understand is the behavior of the
invariant $\SWW$ under compositions of diffeomorphisms.  By reversing
the orientation of all paths, we have (trivially)
\begin{lemma}
           $\SWW(f^{-1},\Gamma) = - \SWW(f,\Gamma)$.
\end{lemma}
The $1$--parameter Donaldson
invariants~\cite{ruberman:isotopy,ruberman:polyisotopy} satisfy the
rule
\begin{equation}\label{hom}
           D(fg) = D(f)+ D(g)
\end{equation}
as does the Seiberg--Witten invariant in the special case that $f$ and
$g$ preserve the $\spinc$ structure.  Each of these statements is
proved by concatenating the relevant paths of metrics or
perturbations, and using the independence of the invariants from the
initial point.  So, for instance, for the Seiberg--Witten invariant,
if $h_0$ is the initial point for calculating $\SW(f,\Gamma)$, then we
use $f^*h_0$ as the initial point for calculating $\SW(g,\Gamma)$.
However, when $f^*\Gamma \not\cong \Gamma$, there seems to be no
choice of paths used for computing $\SWW(f,\Gamma)$, $\SWW(g,\Gamma)$,
and $\SW(fg,\Gamma)$ for which the corresponding moduli spaces can
reasonably be compared.

That the composition rule (if there is one) must be more complicated
than~\eqref{hom} may be seen from the following special case, for
compositions of a diffeomorphism with itself.
\begin{theorem}\label{comp}
           Let $f\co Y\to Y$ be a diffeomorphism such that $\SWW(f,\Gamma)$
           is defined.  Then $\SWW(f_d,\Gamma)$ is defined and satisfies:
	\begin{enumerate}
	\item If $\OO(f,\Gamma)$ is infinite, then $\SWW(f_d,\Gamma) =
           \SWW(f,\Gamma)$ for all $d$.
	\item If $|\OO(f,\Gamma)| = N < \infty$, then:
	$$
	\SWW(f_d,\Gamma) = \frac{\mathrm{lcm}(d,N)}{N} \SWW(f,\Gamma)
	$$
           \end{enumerate}
\end{theorem}
\begin{proof}
           Choose a generic $h_0 \in \Pi$ and path $h$ with which to
           define $\SWW(f,\Gamma)$.  If $|\OO(f,\Gamma)| = \infty$,
           consider the concatenation of paths $H = h \cdot f^*h \cdot
           f_2^*h \cdots f_{d-1}^* h$, with endpoints $h_0$ and $f_d^*
           h_0$.  The union of the moduli spaces $\MM(\Gamma;f_{dk}^*H)$
           corresponding to these paths may be used to define both
           $\SWW(f,\Gamma)$ and $\SWW(f_d,\Gamma)$, and so these
           invariants are equal.

           The same argument works if the orbit of $\Gamma$ is finite,
           except that each moduli space $\MM(\Gamma;f_k^*h)$ is counted
           $\mathrm{lcm}(d,N)/{N}$ times when one takes the union of the
           moduli spaces
           $$
           \MM(\Gamma;f_{dk}^*H)
           $$
           which defines $\SWW(f_d,\Gamma)$.
\end{proof}

We will show later in Corollary~\ref{nonisotopy} that the invariant
$\SWW$ can detect non trivial elements of $\pi_0(\diff_h)$, where
$\diff_h$ is the subgroup of diffeomorphisms which are homotopic to the
identity.  With a better understanding of the invariants of
compositions of diffeomorphisms, it might be possible to recover the
result of~\cite{ruberman:polyisotopy} that $\pi_0(\diff_h)$ is
infinitely generated.
\begin{corollary}\label{isometry}
           Let $f$ be a diffeomorphism which is an isometry for some
           metric on $Y$.  Then $\SWW(f,\Gamma) = 0$ for any $\Gamma$.
\end{corollary}
\begin{proof}  If $f$ were an isometry of a generic metric $g$, then this
           would be clear because one could choose a constant path from
           $g$ to $f^*g$.  Since $g$ might not be generic, we make use
           instead of Theorems~\ref{comp} and~\ref{isotopy}: the isometry
           group of $Y$ is compact, so $f_d$ is in the connected
           component of the identity for some $d$.  But this implies that
           $f_d$ is isotopic to the identity, and so $\SWW(f_d,\Gamma)
           =0$.  By Theorem~\ref{comp}, this implies that $\SWW(f,\Gamma)
           = 0$ as well.
\end{proof}

\section{Calculation of some examples}

In this section we show how to calculate $\SWW$ for the
diffeomorphisms described
in~\cite{ruberman:isotopy,ruberman:polyisotopy}.  These calculations
will be applied to questions about the space $\psc$ in the next
section.  The construction starts with oriented, simply
connected manifolds $X$ satisfying the following conditions.

\begin{itemize}
           \item[(1)] $b^{2}_{+}(X) $ is odd and at least $3$.
           \item[(2)] $\SW(\Gamma_{X}) \neq 0$ for some $\spinc$ structure
           $\Gamma_{X}$ for which $\dim(\MM(\Gamma_{X})) = 0$.
\end{itemize}

The first condition could presumably be weakened to allow
$b^{2}_{+}(X) = 1$; this could be accomplished by developing a version
of $\SWW$ invariants involving a chamber structure.

Given $X$ as above, let $N$ denote the connected sum
$$
\cp^{2} \# \cpbar \# \cpbar
$$
and let $Z = X \# N$.  If $X$ has a homology orientation, then $Z$ 
gets a homology orientation using the natural isomorphism $H^2_+(Z) 
\cong H^2_+(X) \oplus H^2_+(N)$, with the added convention that the 
generator of $ H^2_+(N)$ is dual to the (complex) embedded $2$--sphere 
in $\cp^2$.  In what follows, the homology orientations on $X$ and 
$Z$ will be related in this way, without further mention.

Let $S, E_1,E_2 $ be the obvious embedded $2$--spheres of
self-intersection $\pm 1$ in $N$, and let $\Sigma_\pm = S\pm E_1+E_2$
be spheres of square $-1$.  For $\Sigma = \Sigma_+$ or $\Sigma_-$, let
$\rho^\Sigma$ be the diffeomorphism of $N$ which induces the map
$\rho^\Sigma_*(x) = x + 2(x\cdot \Sigma)\Sigma$ on $H_{2}(N)$, where
`$\cdot$' is the intersection product.  (In cohomology,
$(\rho^\Sigma)^*$ is given by the same formula, where $\Sigma$ is
replaced by $\sigma = PD(\Sigma)$ and $\cdot$ is replaced by $\cup$.)
In particular, letting $s,e_1,e_2$ be the Poincar\'e duals of $S,
E_1,E_2$, we have $(\rho^\Sigma)^*(s) = s + 2(s\cup \sigma)\sigma= 3s
+2(e_2 \pm e_1)$.  Hence $\rho^\Sigma $ is orientation preserving on
$H_+^2(N)$.  Choose the orientation in which $s$ is in the positive
cone.

The diffeomorphism $\rho^{\Sigma}$ glues together with the identity
map of $X$ to give a diffeomorphism $f^\Sigma$ on the connected sum $Z
= X\# N$.  Define $f = f^{\Sigma_+}\circ f^{\Sigma_-}$, and note that
$\alpha(f) = 1$.  Consider the $\spinc$ structure $\Gamma$ on $Z$
whose restriction to $X$ is $\Gamma_{X}$, and whose restriction
$\Gamma_N$ to $N$ has $c_{1} = s+e_{1}+e_{2}$.  It is easy to check
that $\dim(\MM(\Gamma_{Z})) = -1$.
\begin{theorem}\label{swcalc}
           The invariant $\SWW(f,\Gamma_{Z})$ is defined and equals
           $\SW(\Gamma_{X})$.
\end{theorem}
As in~\cite{ruberman:isotopy}, this will be proved by a combination of
gluing and wall-crossing arguments.   As a preliminary to these
arguments, we set up the
notation for wall-crossing.  Consider the hyperbolic space:
$$
\HH = \{x^2 -y^2-z^2 =1, \ x>0\} \subset H^2(N)
$$
For a metric $g^N$, let $\omega({g^N}) \in \HH$ be the unique self-dual
form satisfying
\begin{equation}\label{hyperbolic}
\int_N \omega \wedge \omega = 1 \ \ \mathrm{and} \int_N \omega \wedge
s > 0.
\end{equation}
The \sw\ equations on $N$ admit a reducible solution if and only if
$\omega({g^N})$ lies on the wall $\WW$ in $\HH$
defined by
\begin{equation}\label{walldef}
\int_N \omega \wedge (s+e_{1}+e_{2}) + 2\pi \int_N \omega \wedge \mu = 0.
\end{equation}
The wall is
transversally oriented by the convention that the local intersection
number of a path $\omega_t$ meeting $\WW$ at $t=0$ is given by the
sign of $f'(0)$, where
\begin{equation}\label{wallorient}
f(t) = \int_N \omega_t \wedge (s+e_{1}+e_{2} + 2\pi\mu)
\end{equation}
at $t=0$.

The \sw\ equations may be defined in the context of cylindrical-end
manifolds~\cite{mmr,taubes:l2,morgan-szabo-taubes}, giving rise to
moduli spaces of
finite-energy solutions.  We refer the reader to the
book~\cite{nicolaescu:swbook} for
details; the full strength of the theory is not needed in our case
because the ends of all
the manifolds we consider have positive scalar curvature.   Consider
the manifold
$Y_0$ gotten by removing an open ball from a $4$--manifold $Y$, and
the cylindrical-end
manifold $\hat{Y} = Y_0\, \cup\, S^3
\times [0,\infty)$.   Fix a standard round metric on $S^3$, and
define a  cylindrical-end
metric $\hat{g}$ on $\hat{Y}$ to be one whose restriction to the end
is a product of this
round metric with the usual metric on $[T,\infty)$ for some $T> 0$.
For a pair $\hat{h} =
(\hat{g},\hat{\mu})$, where $\hat{\mu}$ is a self-dual form of
compact support,  one
defines~\cite[\S 4]{nicolaescu:swbook} the moduli space
$\MM(\Gamma^Y,\hat{h})$ consisting
of finite-energy solutions to the \sw\ equations.  Likewise, for a path
$\hat{h}_t$ which is constant in the end, one has a finite-energy
parameterized moduli
space.   In what follows, objects associated to a cylindrical-end
manifold will be
indicated by a $\hat{{}}$ on top.

Because $b_+^2(\hat N) = 1$, the cylindrical-end manifold $\hat{N}$
has wall-crossing
behavior analogous to that on
$N$ itself.   Note the cohomology isomorphism $H^2(\hat{N}) \cong
H^2(N)$, and also note
that classes in $H^2(\hat{N})$ are uniquely represented by $L^2$ harmonic
forms~\cite{aps:I}.  A cylindrical-end metric $\hat{g}$ determines a
unique $L^2$
self-dual harmonic form $\omega({\hat{g}^N}) \in \HH$ satisfying
\eqref{hyperbolic}.  (To
ensure that the integrals in equations \eqref{hyperbolic}, \eqref{walldef} and
\eqref{wallorient} make sense, assume that the generators $e_1,e_2,$
and $s$ are
represented by 2--forms with compact support.)  The wall $\hat{\WW}$
is defined precisely
as above, and the transverse orientation of $\hat{\WW}$ is given by
the sign of the
derivative of the function
\begin{equation}\label{cylwallorient}
\hat{f}(t) = \int_{\hat{N}} \omega_t \wedge (s+e_{1}+e_{2} + 2\pi\hat{\mu})
\end{equation}
as in the compact case.  One readily verifies that there are
reducible finite--energy
solutions to the  \sw\ equations on $\hat{N}$ precisely when
$\omega({\hat{g}^N})\in
\hat{\WW}$.

Beyond the wall-crossing picture, the main ingredient in the
calculation of $\SWW$ is a
special case of the gluing principle for solutions to the
$1$--parameter \sw\ equations.
Cutting and pasting of monopoles (solutions to the usual \sw \
equations) is discussed
comprehensively in section 4.5 of~\cite{nicolaescu:swbook}.  The
analytical details of
gluing in the $1$--parameter case are almost identical.

Write the connected sum $Z = X \# N$ as $X_0 \cup N_0$ where the
union is along the
common $S^3$ boundary.   For any $r >0$, we can view $Z$ as
diffeomorphic to $X_0 \cup
S^3  \times [-r,r]\cup N_0$.  We  consider the perturbed \sw\
equations, where the metric
on $Z$ is required to be glued together from a product metric on the
length $2r$ cylinder
$S^3 \times [-r,r] $ and metrics $g^X \in \Pi(X_0)$ and $g^N \in
\Pi(N_0)$.   Similarly,
we assume that the $g$--self-dual form used in perturbing the \sw\
equations vanishes on
the cylinder.  Such a pair $(g,\mu)$ will be written $h= h^X \#_r
h^N$, and the set of
such pairs will be denoted $\Pi_r(Z)$.   Note that $h \in \Pi_r(Z)$
determines cylindrical
end metrics on $\hat{N}$ and $\hat{X}$, with self-dual $2$--forms of
compact support.
Also, $h \in \Pi_r(Z)$ determines metrics on the closed manifolds $N$
and $X$, by cutting
along $S^3 \times \{0\}$ and gluing in a standard round $4$--ball.

An important remark is that if $r$ is sufficiently large, then the
functions $f$ and
$\hat{f}$ associated to metrics on $N$ determined by $h\in \Pi_r(Z)$
are close in
$\mathcal{C}^1$.  This follows (\cite{clm:I}) from the exponential
decay (cf~\cite{aps:I})
of harmonic forms on $\hat{N}$.  To take advantage of this, consider
a path $h_t \subset
\Pi_r(Z)$ such that $\omega({g_t})$ crosses the wall $\WW$
transversally at $t=0$.  Then
the path
$\omega({\hat{g}_t})$ crosses the wall $\hat{\WW}$ transversally at
$t=0$ and with the
same sign.

The gluing procedure detailed in~\cite[\S 4.5]{nicolaescu:swbook}
relates monopoles on $Z$
to those on $X$ and $N$ in two stages.  The first is to relate
solutions on $Z$ (for
metrics in the space  $\Pi_r(Z)$ where $r$ is sufficiently large) to
finite energy
solutions on $\hat{N}$ and $\hat{X}$.    A second application of the
gluing procedure
relates these, in turn, to solutions on $X$ and $N$, where the
metrics on these manifolds
are as described at the end of the preceding paragraph.

The following basic local calculation leads quickly to the proof of
theorem~\ref{swcalc}.
Note that it evaluates the invariant $\SW(\Gamma;h)$ for the special
case of $h \subset
\Pi_r(Z)$.   Using the invariance properties of $\SWW$, this will
suffice for our
applications.
\begin{proposition}\label{local}  There exists an $r_0$ such that for
all $r > r_0$, the
following statements hold.
           Let $h$ be a path in $\Pi_r(Z)$ such that $h_t = h^X \#_r h^N_t$ for
           all $t$, where $h^X \in \Pi(X)$ is generic and $ h^N_t \subset
           \Pi(N)$ is a generic path.
           \begin{enumerate}
           \item If $\omega({h^N_t})$ is disjoint from $\WW$, then
           $\SW(\Gamma;h) =  0$.
           \item Suppose that $\omega({h^N_t})$ crosses $\WW$ at $t=0$
           transversally, and with positive orientation.  Then
           $\SW(\Gamma;h) = \SW(\Gamma_X)$.
           \end{enumerate}
\end{proposition}

\begin{proof}
To prove these statements,\kern-0.9pt\ we describe the\kern-0.9pt\ $1$--parameter \sw\ moduli space on
$Z$ in terms of the moduli spaces on $\hat N$ and $\hat X$.   The
idea is basically
the same as the `fundamental lemma' in~\cite{fs:swblowup}; details of a
very similar argument
(without the extra parameter) are given in the proof of Theorem 4.5.19
of~\cite{nicolaescu:swbook}.  As a first step, note that the path
$h_t \in \Pi_r(Z)$
induces a path $\hat{h}_t^N$ on
$\hat{N}$ and a constant path
$\hat{h}^X$ on $\hat{X}$.  The assumption that $h^X$ is generic
implies that $\hat{h}^X$
is generic.  Moreover, the moduli space $\MM(\Gamma^X,h^X)$ is
diffeomorphic to the
cylindrical end moduli space associated to $\hat h^X$.

To establish the first part of the proposition, note that if
$\omega({h^N_t})$ is disjoint
from $\WW$, then the same will be true for $\omega({\hat{h}^N_t})$.  Since the
finite-energy
\sw\ moduli space on $\hat{N}$ has formal dimension $-2$ and contains
no reducibles, the
parametrized moduli space
$\MM(\Gamma^{\hat{N}},\hat{h}_t^N)$  will be empty.  Because the
gluing map is a
diffeomorphism, the parameterized moduli space on $Z$ will also be empty.

The second part is a little more complicated; we need to work out the
Kuranishi model for
the $1$--parameter gluing problem for the \sw\ equation for a path of
perturbations for
which the equations on $N$ admit a reducible solution $(A^N_0,0)$ at
$t=0$.   As above,
this path gives rise to a path of perturbations on $\hat{N}$.   Using
the relation between
wall-crossing on $N$ and
$\hat{N}$, there will be exactly one small $\delta$ for which the \sw
\ equations with
perturbation $\hat h_\delta^N$ admit a reducible solution
$(\hat{A}_\delta,0)$.  To simplify
notation, we will assume $\delta = 0$.

Suppose that $(A^X,\varphi^X)$ is a smooth irreducible solution to
the \sw \ equations on
$X$.  As remarked above, it determines a unique finite-energy
solution on $\hat{X}$.
For small $t$, say in an interval $(-\epsilon, \epsilon)$, we have the constant
configuration $(\hat{A}_0,0)$ on $\hat N$. Note that this will be a
solution to the \sw\
equations with perturbation $\hat h^N_t$ only when $t=0$.   Consider
the $1$--parameter
family of configurations on $Z$ gotten by gluing $(A^X,\varphi^X)$ to
$(\hat{A}_0,0)$, and
the problem of deforming this family to give elements in the
$1$--parameter moduli space
on $Z$ parameterized by $h_t = h^X \#_r h^N_t$.   Adapting~\cite[Theorem
4.5.19]{nicolaescu:swbook} to this context, the portion of the
parameterized moduli space
close to the glued-up path near $t=0$ is described (for $r$
sufficiently large) in terms
of the Kuranishi model.  In other words, there is an $S^1$--equivariant map
$$
\psi\co  (-\epsilon,\epsilon) \times \ker(D_{\hat A^N_0}) \times S^1  \to
\coker(D_{\hat A^N_0}) \times H_+^2(\hat N)
$$
such that the parameterized moduli space on $Z$ is diffeomorphic to
$\psi^{-1}(0)/S^1$.
The assumption that the path $\hat h_t$ is generic means that
$\coker(D_{\hat A^N_0})$
vanishes.  From the Atiyah--Patodi--Singer index theorem, the index of
$D_{\hat A^N_0}$ is zero, so that
$$
   \ker(D_{\hat A^N_0}) \cong
\coker(D_{\hat A^N_0}) \cong \{0\}.
$$
As in ~\cite{kronheimer-mrowka:thom}, the derivative $\partial
          \psi/\partial t (0,0,\theta)$ is given by $\hat f'(0)$ where
$\hat f$ is
          defined in equation~\eqref{cylwallorient}.  Hence
$\psi^{-1}(0)/S^1$ consists of
exactly one point; the orientation of this point is the same as the
orientation of the
corresponding point in $\MM(\Gamma^X,h^X)$.   Each point of the moduli space
$\MM(\Gamma^X,h^X)$ contributes therefore $\pm 1$ to $\SW(\Gamma,h)$. 
By part (1), the rest of the path (where $\omega(\hat h_t)$ misses 
$\hat \WW$) doesn't contribute at all, and the result follows.
\end{proof}

Proposition~\ref{local} may be viewed as stating that $\SW(\Gamma;h)$, for
a path which is constant on $X$, is given by the intersection number
of the path $\omega({h^N_t})$ with $\WW$.  Note that this in turn only
depends on the endpoints (as an ordered pair) of the path.
\begin{proof}[Proof of Theorem~\ref{swcalc}]
          Let $h_0 = h^X \# h^N_0 \in \Pi$, where $ h^X$ and $h^N_0$ are
          generic, and choose a generic path $h^N\in \Pi$ from $h^N_0$ to
          $\rho^*(h^N_0)$.  This glues up to give a path $h$ which we can
          use to compute $\SWW(f,\Gamma)$.  According to 
Proposition~\ref{local},
          we need to know the intersection number of the wall $\WW$ with the
          path $(\rho^{\Sigma_+}\circ \rho^{\Sigma_-})_k^*h_N$ for each
          $k$. Now $(\rho^{\Sigma_+}\circ \rho^{\Sigma_-})^* $, viewed as a
          transformation of the hyperbolic space $\HH$, is a parabolic
          element, with fixed point $\imath$ in the unit disc model of
          $\HH$.  In this model, $\WW$ is the geodesic meeting
          the boundary of the unit disc orthogonally at $\imath$ and $1$.
          Hence for any generic starting point $h_0$, there is a unique $n$ for
          which $f_n^*h_0$ is on the left side of $\WW$ and $f_{n+1}^*h_0$
          is on the right side.

          It follows that $\SW(\Gamma,f_k^*h) = 0$ for $k \neq n$, and that
          $\SW(\Gamma,f_n^*h) = \SW(\Gamma_X)$. Hence $\SWW(f,\Gamma) =
          \SW(\Gamma_X)$.
\end{proof}

\section{Isotopy and concordance of positive scalar curvature metrics}

The combination of Theorems~\ref{vanish} and~\ref{comp} give rise to a
method of detecting different components of $\psc$.

\begin{theorem}\label{infinite}
           Suppose that $Y$ is a manifold with a $\spinc$ structure
           $\Gamma$ and diffeomorphism $f\co Y \to Y $ such that
           $\SWW(f,\Gamma) \neq 0$.
           \begin{enumerate}
           \item If $Y$ admits a metric of positive
           scalar curvature, then $\psc(Y)$ has infinitely many components.  In
           fact, if $g_{0}\in \psc$, then the metrics $f_{n}^{*}g_{0}$
           for $n\neq 0$ are all in different path components of $\psc(Y)$.
           \item If $f'$ is another diffeomorphism such that
           $\SWW(f',\Gamma) \neq \SWW(f,\Gamma)$, then for some k, $f_k^*g_0$
           and $(f_k')^*g_0$ are in different path components of $\psc(Y)$.
           \end{enumerate}
\end{theorem}
\begin{proof}
           By Theorem~\ref{comp}, the invariants $\SWW(f_{d},\Gamma)$ are
           all non--zero, as long as $d\neq 0$.  Suppose that
           $f_{l}^{*}g_{0}$ and $f_{k}^{*}g_{0}$ are in the same path
           component of $\psc$.  Hence they are joined by a path $h\co 
           [0,1] \to \psc$.  As in the proof of Theorem~\ref{vanish}, use
           the translates of $h$ by $f_{l-k}$ to show that
           $\SWW(f_{l-k},\Gamma) = 0$.  But this implies that $l=k$.

           For the second part, suppose to the contrary that $f_k^*g_0$
           and $(f_k')^*g_0$ are connected by a path of positive scalar
           curvature metrics for all $k$.  Then by the usual cobordism argument
           we would have
           $$
           \SW(\Gamma; f_{k}^{*}g_{0},f_{k+1}^{*}g_{0}) \neq
           \SW(\Gamma; (f_{k}')^{*}g_{0},(f_{k+1}')^{*}g_{0}).
           $$
           (The assumption that the metrics are connected in $\psc$ provides a
           $2$--parameter family in which the $1$--parameter moduli  spaces
           corresponding to the `sides' are empty.)  This would imply that the
           corresponding $\SWW$ invariants are equal, contradicting our
           assumption.
\end{proof}

\begin{corollary}\label{infinitecomp}{}\qua
	There exist simply-connected $4$--manifolds $Y$ for which\break
	$\psc(Y)$ is non-empty and has infinitely many components.
	More precisely, for any $n \ge 2$ and $k > 10n$, the
	connected sum
	$$
	Z = \#_{2n}\cp^2 \#_k \cpbar
	$$
	has infinitely many components in $\psc(Z)$.
\end{corollary}
\begin{proof}
Let $E$ denote the elliptic surface with $b^2_+ = 2n-1$, blown up at
$k-10n$ points, and let $Z = E \#N$ as in Theorem~\ref{swcalc}.  According
to~\cite{moishezon:sums,mandelbaum-moishezon:algebraic,gompf-stipsicz:book},
$Z$ decomposes as a connected sum as in the statement of the
corollary. Therefore (\cite{gromov-lawson:psc,schoen-yau:psc}) $Z$ admits a
metric $g_0$ of positive scalar curvature.  By Theorem~\ref{swcalc},
$Z$ supports a diffeomorphism with $\SWW(f,\Gamma) = \SW_E(\Gamma)
\neq 0$ (for an appropriate choice of $\Gamma$) and so by
Theorem~\ref{infinite}, $\psc(Z)$ has infinitely many components.
\end{proof}

As mentioned in the introduction, all of the known
obstructions to isotopy are in fact obstructions to the weaker
relation of {\sl concordance}.  By definition, a concordance between
$g_0, g_1 \in \psc(X)$ is a positive scalar curvature metric on $X
\times I$, which in a collar neighborhood of $X \times \{i\}$ is a
product $g_i + dt^2$.  A path in $\psc(X)$, suitably
modified~\cite{gromov-lawson:psc,gajer:cobordism} gives a concordance,
but the converse is an important open question in all dimensions other
than $2$ where it is known to hold.  Although the metrics constructed in
Corollary~\ref{infinitecomp} are not isotopic, it is not at all clear
whether they are concordant.  We present two ways of finding concordant but not
isotopic metrics of positive scalar curvature, corresponding to the
two parts of Theorem~\ref{infinite}. The manifolds in question are the
same as in Corollary~\ref{infinitecomp}, ie, connected sums of
projective space (with both orientations).  The simpler version,
corresponding to the second part of that theorem, uses
the principle that concordant (also known as {\sl pseudoisotopic})
diffeomorphisms give rise to concordant metrics.
\begin{theorem}\label{concordance}
           There are concordant, but not isotopic, metrics of positive
           scalar curvature on simply-connected $4$--manifolds.
\end{theorem}
\begin{proof}
           For any $n \ge 2$, we start with two simply-connected
           $4$--manifolds $X, X'$ satisfying
           \begin{itemize}
	        \item[($1'$)] $b_+^2(X)= 2n-1$.  \item[($2'$)] $X'
           	\simeq X$, inducing a correspondence between $\spinc(
           	X'$ and $\spinc(X)$.  \item[($3'$)] For corresponding
           	$\spinc$ structures $\Gamma_{X'}$ and $\Gamma_{X}$ we
           	have that $\SW(\Gamma_{X'} ) \neq \SW(\Gamma_{X})$.
           	\item[($4'$)] $X' \#\cp^2$ and $X \#\cp^2$ decompose
           	as a connected sum of $\cp^2$'s and $\cpbar$'s.
           \end{itemize}
           There are plenty of sources of such manifolds.  One
           could, for example take $X$ to be a (blown-up) elliptic
           surface, and build the other $X'$ to be the corresponding
           decomposable manifold.  An infinite family $X^{(i)}$ of such
           manifolds arises via the Fintushel--Stern~\cite{fs:knots} knot
           surgery construction on $X$.  The decomposition ($4'$) of the
           stabilized manifold $X^{(i)}\#	\cp^2$ follows readily from the
           analogous fact for $X$.

           Given $X$ and $X'$, we get, as above, diffeomorphisms $f, f'$
           of a single manifold $Z$ with the property that for an
           appropriate $\spinc$ structure on $Z$, the invariants
           $\SWW(f,\Gamma)$ and $\SWW(f',\Gamma)$ are distinct.  This
           $\spinc$ structure has the property that $\OO(f,\Gamma)$ is
           infinite for all $i$.  It is explained
           in~\cite{ruberman:isotopy} and~\cite{ruberman:polyisotopy}
           how to arrange the identification of the stabilized manifolds
           so that $f$ and $f'$ are homotopic to one another.

           Use diffeomorphisms $f,f'$ to pull back a positive scalar
           curvature metric $g_0$ on $Z$ to get metrics $g,g'$.  A
           result of Kreck~\cite{kreck:isotopy} implies that $f$ is
           concordant (or pseudoisotopic) to $f'$.  In other words,
           there is a diffeomorphism $F$ of the $5$--manifold $Z \times
           I$, which is the $f \times id$ in a collar of $Z \times
           \{0\}$ and $f'\times id$ in a collar of $Z\times \{1\}$.
           Pulling back the product metric $g_0 + dt^2$ via $F$ gives
           the required concordance between $g$ and $g'$.  But by the
           second part of Theorem~\ref{infinite}, these metrics are not
           isotopic.
\end{proof}

The nontriviality of the diffeomorphism $f'f_{-1}$ was the main result
of~\cite{ruberman:isotopy}.  Although it is not our main concern in
this paper, we note that the $\SWW$ invariants also recover this
result.
\begin{corollary}\label{nonisotopy}
           The diffeomorphisms $f'_k f_{-k}$ are homotopic to the
           identity, but are not isotopic to the identity.
\end{corollary}
\begin{proof}
           By construction, $f'_k$ is homotopic to $f_k$.  If
           $f'_kf_{-k}$ were isotopic to the identity, then $f'_k$ would
           be isotopic to $f_{k}$.  But, using Theorem~\ref{comp}, this
           would imply that
           $$
           \SWW(f,\Gamma) = \SWW(f_k,\Gamma) = \SWW(f'_k,\Gamma) =
           \SWW(f',\Gamma)
           $$
           which is a contradiction.
\end{proof}

Using some additional topological ingredients, we can prove a stronger
result highlighting the difference between concordance and isotopy of
$\psc$ metrics.
\begin{theorem}\label{concisotopy}
          Let $Z$ be a connected sum $\#_{2n}\cp^2 \#_k \cpbar$ as in
          theorem~\ref{infinitecomp}.  Then any concordance class of
          $\psc$ metrics on $Z$ contains infinitely many isotopy
          classes.
\end{theorem}
\begin{proof}
          The main topological ingredient is the computation of the
          group $\Delta_n$ of bordisms of diffeomorphisms.  Recall that
          $\Delta_n$ consists of pairs $(X,f)$ where $f$ is a
          diffeomorphism of $X$, modulo cobordisms over which the
          diffeomorphism extends.  (All manifolds are oriented, and
          cobordisms and diffeomorphisms are to respect the
          orientations.)  The group was computed by
          Kreck~\cite{kreck:bordism} and Quinn~\cite{quinn:bordism}; we
          will recall Kreck's computation of $\Delta_4$ shortly.  To
          make use of this computation, we need the following
          observation.
\begin{claim}
          Let $A$ be a non-spin, simply-connected $4$--manifold, with a
          metric $g_0$ of positive scalar curvature.  Suppose that $f$
          is a diffeomorphism which is bordant to $id_A$.  Then $f^*g_0$
          is concordant to $g_0$.
\end{claim}
\begin{proof}[Proof of Claim]
          By definition, there is a manifold $B^5$ with two boundary
          components $\partial_0B \cong A \cong \partial_1B$, and a
          diffeomorphism $F\co B \to B$ with
          $$
          F|\partial_0B = id_A \quad \mathrm{and}\quad F|\partial_1B =
          f.
          $$
          According to~\cite[Remark 11.5]{kreck:bordism} we can (and
          will) assume that $B$ is simply connected.  Now there is a
          cobordism $W^6$, relative to the boundary of $B$, to the
          product $A \times I$.  Again, by preliminary surgery on
          circles if necessary, we can assume that $W$ is simply
          connected. By surgery on $2$--spheres in $W$, we can also
          assume that
          $$
          \pi_2(W,B) = \pi_2(W,A \times I) = 0.
          $$
          (It is at this point that we used the assumption that $A$ is
          not spin, in order to be able to choose a basis of
          $2$--spheres with trivial normal bundle.)  As a consequence,
          we can find a handle decomposition for $W$, relative to the
          boundary, with only $3$--handles.  In other words, $B$ is
          obtained from $A \times I$ by surgery on $2$--spheres, and
          {\sl vice-versa}.

          Now $A \times I$ has the product metric $g_0 + dt^2$, which
          can be pushed across the cobordism $W$ to give a metric of
          positive scalar curvature on $B$ which is a product near the
          boundaries.  Use the diffeomorphism $F$ to pull this metric
          back to give a new positive scalar curvature metric on $B$
          which is $g_0 + dt^2$ near $\partial_0 B$ and $f^* g_0 + dt^2$
          near $\partial_1 B$.  Push this new metric {\sl back} across
          $W$ to give a positive scalar curvature metric on $A \times
          I$; since the metric is not changed near the boundary
          components, this is a concordance.
          (Cf~\cite{gajer:concordance} for an argument of this sort.)
          \renewcommand{\qedsymbol}{$\sqr55_{\;\mathrm{Claim}}$}
\end{proof}
\def\qedsymbol{$\sqr55$}

          To make use of the claim, we need to show that the
          diffeomorphism $f$ on $Z = Z \# N \cong \#_{2n}\cp^2
          \#_k\cpbar$ constructed in Theorem~\ref{swcalc} is bordant to
          the identity.  By the basic theorem of~\cite{kreck:bordism},
          this is determined by the action of $f_*$ on $H_2(Z)$, as
          follows. Let $q$ denote the intersection form on $H_2(Z)$.
          Then $f$ is bordant to $id_Z$ if (and only if) there is a
          $1/2$--dimensional subspace
          $$
          V \subset H_2(Z) \oplus H_2(Z)
          $$
          which is $f_* \oplus id$ invariant and on which the form $q
          \oplus -q$ vanishes. ($V$ is called a metabolizer for the
          isometric structure $(H_2(Z) \oplus H_2(Z),f_* \oplus
          id,q\oplus -q).$) This will clearly hold if there is a
          metabolizer for the action of $f_* \oplus id $ on the summand
          $H_2(N) \oplus H_2(N)$ of $H_2(Z) \oplus H_2(Z)$.  Now a
          straightforward calculation shows that on $H_2(N)$, with
          respect to the basis $\{S,E_1,E_2\}$ the action of $f_*$ is
          given by the matrix:
          $$
\left (\begin {array}{ccc} 9&4&-8\\\noalign{\medskip}4&1&-4
\\\noalign{\medskip}8&4&-7\end {array}\right)
          $$
          The reader can check that with respect to the basis given above,
the vectors
$(1,0,1,0,0,0),
          (0,1,0,0,1,0),$ and $(1,0,1,1,0,1)$
           span a metabolizer in $H_2(N) \oplus H_2(N)$.  It
          follows that $f$ is bordant to the identity.

          To conclude the proof, let $g_0$ be a representative of any
          concordance class in $\psc(Z)$.  (Because of the
          decomposition of $Z$ as a connected sum, this space is
          non-empty.)  Then all of the metrics $f_k^*g_0$ are
          concordant, by the claim and the preceding calculation.  On
          the other hand, these metrics are mutually non-isotopic, because of
          Theorem~\ref{infinite}.
\end{proof}
\section{The moduli space of positive scalar curvature metrics}
As mentioned in the introduction, one can also study the space of
metrics modulo the diffeomorphism group.  From this point of view, the
different path components of $\psc$ detected in the previous sections
are really the same because they are related by a diffeomorphism.
\begin{defn}
          Let $X$ be a manifold for which $\psc(X)$ is nonempty.  Then we
          define the moduli space $\mpsc(X)$ to be $\psc(X)/\diff(X)$, where
          the action of $\diff(X)$ is by pullback.
\end{defn}
The main goal of this section is to show that $\mpsc(X)$
can have some nontrivial topology.
\begin{theorem}\label{disconnect}
          For any $N$, there is a smooth $4$ manifold $X$ such that
          $\mpsc(X)$ has at least $N$ components.
\end{theorem}
The manifolds we construct to prove Theorem~\ref{disconnect} are
non-orientable.  It would be of some interest to find orientable, or
even simply connected $4$--manifolds with $\mpsc$ disconnected; these
would have no known analogues among higher (even) dimensional
manifolds.  At the end of the paper, we will suggest some potential examples.

\subsection{Construction of examples}
The phenomenon underlying Theorem~\ref{disconnect} is that the same
$C^\infty$ manifold can be constructed in a number of ways.  Each
construction gives rise to a metric of positive scalar curvature;
eventually we will show that these metrics live in different
components of $\mpsc$.  The basic ingredient is a non-orientable
manifold described in~\cite{gilkey:even}.  Start with the
$3$--dimensional lens space $L(p,q)$, which carries a metric of
positive curvature as a quotient of $S^3$.  Assuming, without loss of
generality, that $q$ is odd, there is a double covering
$$
L(p,q) \to L(2p,q)
$$
with covering involution $\tau\co L(p,q) \to L(p,q)$.  Let
$$
M(p,q) = \left( S^1 \times L(p,q) \right)/((z,x) \sim (\bar{z},\tau(x))).
$$
This is a non-orientable manifold, because complex conjugation on the
$S^1$ factor reverses orientation.  The covering translation is an
isometry of the product metric on  $S^1 \times L(p,q)$, and so descends
to a metric of positive scalar curvature on  $M(p,q)$.  A second
useful description of $M(p,q)$ results from the projection of $S^1
\times L(p,q)$ onto the second factor.  The descends to an $S^1 $
fibration of $M(p,q)$ over $L(2p,q)$.  From the fibration over $S^1$,
the fundamental group of $M(p,q)$ is given as an extension
$$
       \{1\}  \to \Z/p \times \Z \to \pi_1(M(p,q)) \to \Z/2 \to \{1\}
$$
whereas the second description yields the exact sequence:
\begin{equation}\label{pi1M}
       \{1\}  \to \Z \to \pi_1(M(p,q)) \to \Z/(2p) \to \{1\}
\end{equation}
\begin{proposition}
          If $M(p,q) \cong M(p',q')$, then $L(p',q') \cong L(p,q)$.
\end{proposition}
This is straightforward; if the manifolds are diffeomorphic, then
their orientable double covers are diffeomorphic.  But it is a
standard fact that $S^1 \times L(p,q)$ determines the lens space up to
diffeomorphism.

In view of this proposition, it is somewhat surprising that most of
the information about $L$ (except for its fundamental group) is lost
when one does a single surgery on $M$.

Let $\gamma \subset M$ be a fiber over a point of $L(2p,q)$; by
construction it is an orientation preserving curve.  Note that we can
arrange that $\gamma$ meets every lens space fibre in the fibration
$M(p,q) \to S^1$ transversally in one point.  Choose a trivialization
of the normal bundle of $\gamma$, and do surgery on it, to obtain a
manifold $X(p,q)$.  The analogous oriented construction (ie, surgery
on $S^1 \times \Sigma^3$) is called the spin of $\Sigma$, so we will
call $X(p,q)$ the {\em flip-spun} lens space.  From the exact
sequence~\eqref{pi1M}, we see that $\pi_1(X) \cong \Z/(2p)$.
Moreover, the evident cobordism $W(p,q)$ between $M(p,q)$ and $X(p,q)$
has fundamental group $\Z/(2p)$.  The inclusion of $X(p,q)$ into this
cobordism induces an isomorphism on $\pi_1$, while the inclusion of
$M(p,q)$ induces the projection in the sequence~\eqref{pi1M}.

By studying the detailed properties of torus actions on these
manifolds, P~Pao \cite{pao:torus-I} showed that the spin of the lens
space $L(p,q)$ is (non-equivariantly) diffeomorphic to the spin of
$L(p,q')$.  We show that the same holds true for the flip-spun lens
spaces. The proof is an adaptation of an argument for the orientable
case; the $4$--dimensional version was explained to me by S~Akbulut,
and the $5$--dimensional version was known to I~Hambleton and
M~Kreck.
\begin{theorem}\label{flip}
          The manifolds $X(p,q)$ and $X(p,q')$ are diffeomorphic for any
          choice of $q$ and $q'$, and are independent of the choice of
          framing of the surgery circle.
\end{theorem}
As a consequence of the theorem, we will refer to all of these
manifolds as $X(p)$.  Note that the diffeomorphism described in the
proof below respects of the generator of $\pi_1(X(p,q))$ corresponding
to the core of the $1$--handle.  So there is also a preferred
generator of $\pi_1(X(p))$.
\begin{proof}
          There are two framings for the surgery on $\gamma$, which differ
          by twisting by the non-trivial element in $\pi_1(\SO(3))$.  But
          since $L(p,q)$ has a circle action with a circle's worth of
          fixed points, this rotation can be undone via an isotopy of
          $M(p,q)$ with support in a neighborhood of a copy of $L(p,q)$
          of $M(p,q)$.  Thus the two surgeries give the same manifold.

          The key observation for showing that $X(p,q)$ does not depend on
          $q$ is that it has a very simple structure as a double of a
          manifold with boundary.  Here is a more precise description:
          First note that splitting the $S^1$ fibers into pieces with
          $\mathrm{Re}(z) \ge 0$ or $\mathrm{Re}(z) \le 0$ exhibits
          $M(p,q)$ as the double of the non-orientable $\I$--bundle
          $L(2p,q) \ttimes \I$.  $L_0$ be $L(2p,q)$ minus an open disc,
          and consider the restriction $L_0 \ttimes \I$ of this
          $\I$--bundle.  Then $X(p,q)$ is the double of $L_0 \ttimes \I$.

          To see this, note that the curve $\gamma$ along which the surgery
          is done is the union of two $\I$ fibres whose endpoints are
          identified when $L(2p,q) \ttimes \I$ is doubled.  $M(p,q)$ minus
          a neighborhood of $\gamma$ is clearly gotten by doubling $L_0
          \ttimes \I$ along the part of its boundary lying in the boundary
          of $L(2p,q) \ttimes \I$.  When doing the surgery, one glues in a
          copy of $D^2$ for each circle in the boundary of a neighborhood
          of $\gamma$.  This has the effect of the doubling along the rest
          of the boundary of $L_0 \ttimes \I$.

          The rest of the argument is handlebody theory, and can be
          expressed in either $4$ or $5$--dimensional terms. Both versions
          start with the observation that $L_0$ has a handlebody
          decomposition with a $0$--handle, a $1$--handle, and a single
          $2$--handle.  $X(p,q)$, being the double of $L_0 \ttimes \I$, is
          the boundary of the $5$--manifold $L_0 \ttimes \I \times \I$,
          which has a handlebody decomposition with handles which are
          thickenings of those of $L_0$ and therefore have the same
          indices.  But the isotopy class of the attaching map of the
          $2$--handle, being an embedding of a circle in the $4$--manifold
          $S^1 \ttimes S^3$, depends only on its homotopy class and its
          framing.  By considering the $4$--dimensional picture, the
          framing is also seen to be independent of $q$.  But this
          homotopy class is $2p$ times a generator, and is therefore
          independent of $q$. By considering the $4$--dimensional picture,
          the framing is also seen to be independent of $q$.  It follows
          that $X(p,q) \cong X(p,q')$ for any $q,q'$.

          The $4$--dimensional version is a little more subtle; the fact
          that $X(p,q)$ is a double implies that it has a
          ($4$--dimensional) handle decomposition with a (non-orientable)
          $1$--handle, two $2$--handles, a $3$--handle, and a $4$--handle.
          As usual, we can ignore the handles of index $> 2$.  The $0$ and
          $1$--handles, plus the first $2$--handle, give a handlebody
          picture of $L_0 \ttimes \I$ as in~\cite[Figure
          4.39]{gompf-stipsicz:book}.  If one draws this first $2$--handle
          as in Figure 2 below, then it has framing $2pq$, because the
          framing is the same as the normal framing of the $(2p,q)$ torus
          knot which is the attaching region in the usual Heegaard
          splitting of $L(2p,q)$.  The $1$--handle is of course
          non-orientable; this means that in the diagram, a curve going
          into the front of one of the attaching spheres of the
          $1$--handle comes out the back of the other one.  The second
          $2$--handle is attached, with framing $0$, along a meridian to
          the first $2$--handle.  By sliding the first $2$--handle over the
          second, one can change crossings; however this is not quite
          enough to change the diagram for $X(p,q)$ into the one for
          $X(p,q')$.

          The reason for this is that in the diagram for $X(p,q)$ the
          $2$--handle is drawn as a tangle of $p$ arcs going from one
          attaching sphere of the $1$--handle to the other, and changing
          crossings does not change the endpoints of this
          tangle. Associated to the tangle is a permutation of the $2p$
          points where the $2$--handle crosses the attaching sphere of the
          $1$--handle.  For differing $q,q'$, these may not be the same
          permutation.  However they are both $2p$--cycles and are
          therefore conjugate.  Choose a permutation conjugating one to the
          other, and a $2p$--string braid $\sigma$ inducing this
          permutation on the endpoints of the strings.  Without changing
          $X(p,q)$, we can alter the $2$--handle by inserting $\sigma$ on
          one end and $\bar{\sigma}$ on the other end, where $\bar{\sigma}$
          is $\sigma^{-1}$ with all of the crossings reversed. (Note that
          $\sigma^{-1}$ becomes $\bar{\sigma}$ when pushed over the
          $1$---handle.)  Now we can use the second $2$---handle to change
          crossings until the diagrams for the two manifolds are
          identical.  Since $2pq$ and $2pq'$ are both even, we can arrange
          by further handle slides to make the framings the same.
\end{proof}

The diagram below shows the proof for the diffeomorphism $X(4,1) \cong
X(4,3)$.  The top picture is a handle decomposition for $X(4,1)$.  The
second is equivalent to the first (push the braid at the left over the
$1$--handle and cancel with the braid on the right).  The result
(exercise!) is the diagram for $X(4,3)$, except that the framing is
$8$ rather than $24$.  But, as in the proof, this can be remedied by
sliding $8$ times over the second $2$--handle.  Because $2$ is a small
value for $p$, one does not have to change any crossings.

\begin{figure}[ht!]
\begin{center}
          \vskip2ex
          \includegraphics{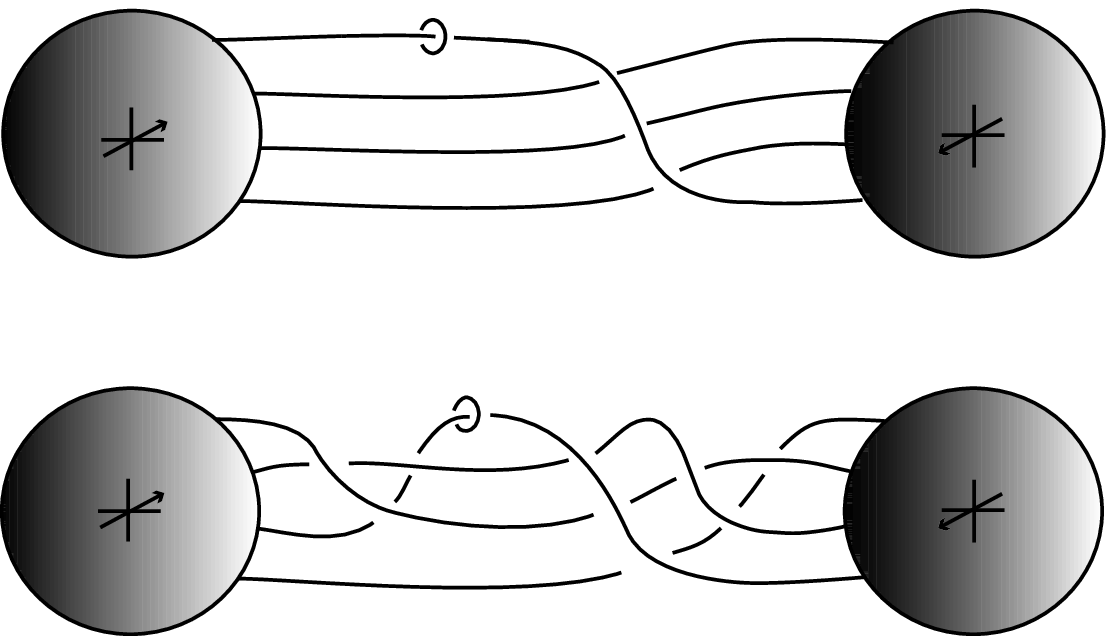}
          \\\small Figure 2
\end{center}
\end{figure}

Corresponding to each of its descriptions as surgery on $M(p,q)$, the
manifold $X(p)$ has a metric of positive scalar curvature.  As
described above, the manifold $M(p,q)$ has a positive scalar curvature
metric $g_{p,q}$ which is locally a product of a metric on $S^1$ with the
positive curvature metric on $L(p,q)$ which it inherits from $S^3$.
According to the surgery construction of
Gromov-Lawson/Schoen-Yau~\cite{gromov-lawson:psc,schoen-yau:psc}, if
the normal disc bundle to $\gamma$ has sufficiently small radius, then
$X(p)$ acquires a Riemannian metric $g_{p,q}$ of positive scalar
curvature.  (Technically speaking, this metric is only well-defined
up to certain choices; none of these will affect what follows.)  Note
that $g_{p,q}$ depends, at least in principle, on $q$ and not just on
$p$.  Moreover~\cite{gajer:cobordism}, the cobordism $W(p,q)$ between
$M(p,q)$ and $X(p)$ has a metric of positive scalar curvature
which is a product along both boundary components.

\subsection{$\eta$--invariants for twisted $\mathbf{\pinpm}$ structures}
We briefly review the invariant we will use to detect the different
components of $\mpsc$.  Let $Y$ be a Riemannian manifold with a
$\pinpm$ structure. A unitary representation $\alpha$ of $\pi_{1}(Y)$
gives rise to a twisted Dirac $D$ operator on $Y$.  Extending the work
of Atiyah--Patodi--Singer~\cite{aps:I} to this context,
Gilkey~\cite{gilkey:even} proved that the $\eta$ function associated
to $D$ had an analytic continuation which was regular at the origin,
and defined the invariant
$$
\eta(Y,g,\alpha) = \frac{1}{2}(\eta(0,D) + \dim\ker(D)).
$$
(We have omitted any notation specifying the $\pinpm$ structure.)

When $Y$ is even-dimensional, $\eta(Y,g,\alpha)$ modulo $\Z$ is
independent of the metric, and thus is an $\R/\Z$--valued invariant
of the underlying smooth manifold (even up to an appropriate notion of
bordism).  On the other hand, as a function of the metric, it is
well-defined in $\R$.  From the Atiyah--Patodi--Singer~\cite{aps:I}
theorem, plus the Lichnerowicz formula, one deduces a stronger
invariance property.
\begin{lemma}[\cite{botvinnik-gilkey:spaceforms}, Lemma 2.1]\label{extend}
        Let $W$ be a $\pinpm$ manifold with boundary $Y \coprod Y'$, and
        suppose that there is a unitary representation of $\pi_1(W)$
        restricting to representations $\alpha,\alpha'$ on the boundary. If
        $W$ has a metric of positive scalar curvature which is a product near
        its boundary, then
        $$
        \eta(Y,g,\alpha) = \eta(Y',g',\alpha').
        $$
\end{lemma}
In particular, metrics of positive scalar curvature which are in the
same path component of $\psc(Y)$ have the same $\eta$--invariant, for
every unitary representation of $\pi_1(Y)$.
\subsection{Proof of Theorem~\ref{disconnect}}
We will actually prove the following more precise version of
Theorem~\ref{disconnect}.
\begin{theorem}\label{lensdisconnect}
Suppose that $p$ is odd.  If the metrics $g_{p,q}$ and $g_{p,q'}$ on
$X(p)$ lie in the same component of $\mpsc(X(p))$, then $q' \equiv
(q)^{\pm 1} \pmod{2p}$.
\end{theorem}
By taking $p$ sufficiently large, so that there are many residue
classes in $\Z_{2p}^*$, we get arbitrarily many components in
$\mpsc(X(p))$ as stated in Theorem~\ref{disconnect}.  Note that the
converse of Theorem~\ref{lensdisconnect} is trivially true.
\begin{proof}[Proof of Theorem~\ref{lensdisconnect}]
According to Lemma~\ref{extend}, if $g_{p,q}$ and $g_{p,q'}$ are in
the same component of $\mpsc(X(p))$, the corresponding
$\eta$--invariants must coincide.  We will show that this implies the
stated equality of $q'$ and $q^{\pm 1}$.  First we must compute the
$\eta$--invariants.
\begin{claim}
        For at least one choice of framing, the cobordism $W$ is a $\pinp$
        cobordism from $M(p,q)$ to $X(p)$.
\end{claim}
\begin{proof}[Proof of claim]
        The obstruction to the extending a $\pinp$ structure over the
        cobordism is a relative $w_2$, lying in $H^2(W,M;\Z/2)$.  This is the
        same as the obstruction to extending the restriction of the $\pinp$
        structure to a neighborhood $N$ of the attaching curve $\gamma$ over
        the $2$--handle.  There are two $\pinp$ structures on $N$, which differ
        by twisting by the non-trivial element of $\pi_{1}\SO(3)$; one of them
        extends over the $2$--handle.  If, with respect to one framing,
        the $\pinp$ structure does not extend, then it will if the handle
        is attached by the other framing, which also differs by twisting by
        the non-trivial element of $\pi_{1}\SO(3)$.
\end{proof}

Fix a $\pinp$ structure on $M(p,q)$, and choose the framing of the
$2$--handle so that $W$ gives a $\pinp$ cobordism to $X(p)$.  Let
$\alpha_s$, $s = 0,\ldots,2p-1$ be the $U(1)$
representations of $\pi_1(X(p)$.  By extending over the cobordism
$W(p,q)$, these give rise to $U(1)$ representations of
$\pi_1(M(p,q))$, for which we will use the same name.  Each of these
representations gives rise to a twisted $\pinp$ structure on $W(p,q)$.
According to Lemma~\ref{extend}:
\begin{equation}\label{equaleta}
        \eta(X(p),g_{p,q},\alpha_s) =\eta(M(p,q),g_{p,q},\alpha_s)
\end{equation}

A  $\pinp$ structure twisted by a $U(1)$ representation is in a
natural way a $\pinc$ structure, and so the $\eta$ invariants
in~\eqref{equaleta} may be calculated as in~\cite[Theorem
5.3]{gilkey:even}:
\begin{equation}\label{pinceta}
            \eta(X(p),g_{p,q},\alpha_s)  =
            \rho_{\alpha_s}(L(2p,q)) -
            \rho_{\alpha_{s+p}}(L(2p,q))
\end{equation}
where the reduced $\eta$--invariant $ \rho_{\alpha_s}(L(2p,q))$ (in the
original notation of~\cite{aps:II}, but not of~\cite{gilkey:even}) is given by:
$$
        \frac{1}{2p}
        \sum_{\stackrel{\lambda^{2p} = 1}{\lambda \neq 1}}
        \frac{(\lambda^s-1)(\lambda^q)}{(\lambda^q-1)(\lambda -1)}
$$
Note that with the given metrics, these formulas are valid over $\R$,
not just reduced modulo $\Z$.  This follows, for example, from
\cite{gilkey:sphere} where $ \rho_{\alpha_s}(L(2p,q))$ is computed via
a spectral decomposition for the relevant Dirac operators in terms of
spherical harmonics.  Thus:
\begin{equation}
        \eta(X(p),g_{p,q},\alpha_s)  =
        \frac{1}{2p} \sum_{\stackrel{\lambda^{2p} = 1}{\lambda \neq 1}}
        \frac{(\lambda^s-\lambda^{s+p})(\lambda^q)}{(\lambda^q-1)(\lambda -1)}
\end{equation}
Dividing the $2p^{th}$ roots of unity into those for which $\lambda^p
= \pm 1$, this may be rewritten as:
$$
       \eta(X(p),g_{p,q},\alpha_s)  =
            \frac{1}{p} \sum_{\lambda^{p} = -1}
            \frac{(\lambda^s)(\lambda^q)}{(\lambda^q-1)(\lambda -1)}
$$
When $p$ is odd, we can substitute $-\lambda$ for $\lambda$ to get
$$
\eta(X(p),g_{p,q},\alpha_s)
            = \frac{1}{p}
            \sum_{\lambda^{p} = 1}
            \frac{(-1)^{s+1}(\lambda^s)(\lambda^q)}{(\lambda^q+1)(\lambda +1)}
$$
where we have also used the fact that $q$ is odd.

Suppose now that $g_{p,q}$ and $g_{p,q'}$ are connected in
$\mpsc(X(p))$.  Taking into account that a diffeomorphism might
permute the $\pinc$ structures on $X(p)$, this means that for some $a
\in \Z/(2p)^*$, we must have
\begin{equation}\label{permeta}
       \eta(X(p),g_{p,q},\alpha_s) =  \eta(X(p),g_{p,q'},\alpha_{as})\quad
       \forall s = 0,\dots,2p-1.
\end{equation}
Under the assumption that $p$ is odd, we will show that
equation~\eqref{permeta} implies that either $q \equiv q' \pmod{2p}$
with $a=1$, or $qq' \equiv 1 \pmod{2p}$ with $a = -q'$.  Since $p$ is
odd, and we have assumed that $q,q'$ are odd as well, it suffices to
show that one of these congruences holds mod $p$.

Let $\omega = e^{2 \pi i/p}$ be a primitive $p^{th}$ root of unity,
and note that for any $j= 1,\dots,p-1$ we have
\begin{align*}
\sum_{s=0}^{p-1} (-1)^{s+1}\eta(X(p),g_{p,q},\alpha_s) \omega^{-js}
    &= \frac{1}{p} \sum_{s=0}^{p-1} \sum_{k=0}^{p-1}
\frac{\omega^{(k-j)s}\omega^{kq}}{(\omega^{kq} + 1)(\omega^{k} + 1)}\\
&=
\frac{1}{p} \sum_{k=0}^{p-1}
\frac{\omega^{kq}}{(\omega^{kq} + 1)(\omega^{k} + 1)}
\sum_{s=0}^{p-1} \omega^{(k-j)s}\\
&= \frac{\omega^{jq}}{(\omega^{jq} + 1)(\omega^{j} + 1)} =
\frac{1}{(\omega^{-jq} + 1)(\omega^{j} + 1)}\\
&=
\frac{(\omega^{-2jq}-1)(\omega^{2j}-1)}{(\omega^{-jq}-1)(\omega^{j}-1)}
\end{align*}

Suppose that the metrics $g_{p,q}$ and $g_{p,q'}$ are connected in
$\mpsc(X)$, so that the $\eta$--invariants correspond as in
equation~\eqref{permeta}.  By repeating these same calculations for $
\eta(X(p),g_{p,q'},a\alpha_s)$, we get that
$$
\frac{(\omega^{-2jq}-1)(\omega^{2j}-1)}{(\omega^{-jq}-1)(\omega^{j}-1)}
=
\frac{(\omega^{-2bjq'}-1)(\omega^{2bj}-1)}{(\omega^{-bjq'}-1)(\omega^{bj}-1)}
$$
where $ab \equiv 1 \pmod{p}$.

We now use Franz's independence lemma~\cite{atiyah-bott:lefschetz} which
states that there are no non-trivial multiplicative relations amongst
the algebraic numbers $\omega^k - 1$, ie, the terms on either side of
this equality must match up in pairs.  Eliminating some trivial
possibilities we find that either $b = -q$ (implying $qq' \equiv 1
\pmod{p}$) or $b=1$ which implies $q \equiv q' \pmod{p}$.
Theorem~\ref{lensdisconnect} follows.
\end{proof}
\section{Concluding Remarks}
The arguments in the previous section mostly extend to twisted
products of higher-dimensional lens spaces; and so one might be able
to use our construction to get explicit higher-dimensional manifolds
with $\mpsc$ disconnected.  The only missing step is the
diffeomorphism between the different flip-spun lens spaces.  The
handlebody arguments do not work in higher dimensions, so one would
have to resort to surgery theory.

As mentioned earlier, the spin of $L(p,q)$ ($S(L)$, given by surgery
on $S^1 \times L$) is diffeomorphic to the spin of $L(p,q')$.  Hence
$S(L)$ supports a variety of metrics $g_{p,q}$ of positive scalar
curvature, and it is natural to speculate that these are in different
components of $\mpsc$.  Unfortunately, there does not seem to be a
non-trivial $\eta$--invariant for even-dimensional orientable
manifolds which will detect these components.  In principle, the
Seiberg--Witten equations could give rise to obstructions to homotopy
of these metrics in $\psc$ or $\mpsc$, because the dimension of the
\sw\ moduli space (for any $\spinc$ structure on $S(L)$) is $-1$.  The
definition of $1$--parameter invariants, following the scheme in the
first part of this paper, becomes complicated by the presence of a
reducible solution.  We hope to resolve these issues in future work.


\begin{thebibliography}

\bibitem{atiyah-bott:lefschetz} {\bf M\,F Atiyah}, {\bf R~Bott},
\emph{A {Lefschetz} fixed point formula for elliptic complexes: {II}},
Annals of Math. {88} (1968) 451--491

\bibitem{aps:I} {\bf M\,F Atiyah}, {\bf V\,K Patodi}, {\bf
I\,M Singer}, \emph{Spectral asymmetry and {Riemannian} geometry:
{I}}, Math.\ Proc.\ Camb.\ Phil.\ Soc. {77} (1975) 43--69

\bibitem{aps:II} {\bf M\,F Atiyah}, {\bf V\,K Patodi}, {\bf
I\,M Singer}, \emph{Spectral asymmetry and {Riemannian} geometry:
{II}}, Math.\ Proc.\ Camb.\ Phil.\ Soc. {78} (1975) 405--432

\bibitem{botvinnik-gilkey:spaceforms} {\bf B Botvinnik}, {\bf P\,B
Gilkey}, \emph{Metrics of positive scalar curvature on spherical space
forms}, Canad. J. Math. {48} (1996) 64--80

\bibitem{clm:I} {\bf S Cappell}, {\bf R Lee}, {\bf E~Miller},
\emph{Self-adjoint elliptic operators and manifold
decompositions. {Part I}: {Low} eigenmodes and stretching}, Comm.\
Pure Appl.\ Math. {49} (1996) 825--866

\bibitem{fs:swblowup}
{\bf R Fintushel}, {\bf R\,J Stern},  \emph{Immersed spheres in
$4$--manifolds and the immersed
Thom conjecture}, Turkish J. Math. {19} (1995) 145--157

\bibitem{fs:knots} {\bf R Fintushel}, {\bf R\,J Stern}, \emph{Knots,
links and $4$--manifolds}, Invent. Math. {134} (1998) 363--400

\bibitem{gajer:cobordism} {\bf P Gajer}, \emph{Riemannian metrics of
positive scalar curvature on compact manifolds with boundary},
Ann. Global Anal. Geom.  {5} (1987) 179--191

\bibitem{gajer:concordance} {\bf P Gajer}, \emph{Concordances of
 metrics of positive scalar curvature}, Pacific J. Math. {157} (1993)
 257--268

\bibitem{gilkey:even} {\bf P\,B Gilkey}, \emph{The eta invariant for
even-dimensional {${\rm PIN}\sb {{\rm c}}$} manifolds}, Adv. in
Math. {58} (1985) 243--284

\bibitem{gilkey:sphere} {\bf P\,B Gilkey}, \emph{The geometry of
spherical space form groups}, World Scientific Publishing Co. Inc.
Teaneck, NJ (1989) with an appendix by A Bahri and M Bendersky

\bibitem{gompf-stipsicz:book} {\bf R\,E Gompf}, {\bf A\,I Stipsicz},
\emph{$4$--manifolds and {K}irby calculus}, American Mathematical
Society, Providence, RI (1999)

\bibitem{gromov-lawson:psc} {\bf M Gromov}, {\bf H\,B Lawson Jr},
\emph{The classification of simply connected manifolds of positive
scalar curvature}, Ann. of Math. {111} (1980) 423--434

\bibitem{kreck:isotopy} {\bf M Kreck}, \emph{Isotopy classes of
diffeomorphisms of $(k-1)$--connected almost-parallelizable
$2k$--manifolds}, from: ``Algebraic topology, Aarhus 1978 (Proc.  Sympos.
Univ. Aarhus, 1978)'', Springer, Berlin (1979) 643--663

\bibitem{kreck:bordism} {\bf M Kreck}, \emph{Bordism of
diffeomorphisms and related topics}, Springer--Verlag, Berlin (1984)
with an appendix by Neal W Stoltzfus

\bibitem{kronheimer-mrowka:thom} {\bf P\,B Kronheimer}, {\bf T\,S
Mrowka}, \emph{The genus of embedded surfaces in the projective
plane}, Math.\ Res.\ Lett. {1} (1994) 797--808

\bibitem{lawson-michelson} {\bf H\,B Lawson Jr}, {\bf M-L Michelsohn},
\emph{Spin geometry}, Princeton University Press, Princeton, NJ (1989)

\bibitem{lichnerowicz:spinors} {\bf A Lichnerowicz}, \emph{Spineurs
harmoniques}, C. R. Acad. Sci. Paris {257} (1963) 7--9

\bibitem{mandelbaum-moishezon:algebraic}
{\bf R Mandelbaum}, {\bf B Moishezon}, \emph{On the topology of simply
     connected algebraic surfaces}, Trans. Amer. Math. Soc. {260} (1980)
     no.~1, 195--222

\bibitem{mmr} {\bf J~Morgan}, {\bf T~Mrowka}, {\bf D~Ruberman},
\emph{The ${L}^2$--moduli space and a vanishing theorem for
{Donaldson} invariants}, Monographs in Geometry and Topology {2},
International Press (1994)

\bibitem{moishezon:sums} {\bf B Moishezon}, \emph{Complex surfaces and
connected sums of complex projective planes}, Lecture Notes in
Mathematics, 603, Springer--Verlag, Berlin (1977) with an appendix by
R Livne

\bibitem{morgan-szabo-taubes} {\bf J~Morgan}, {\bf Z~Szabo}, {\bf
C~Taubes}, \emph{A product formula for the Seiberg--Witten invariants
and the generalized Thom conjecture}, J.  Diff.\ Geo. {44} (1996)
706--788

\bibitem{nicolaescu:swbook}
{\bf L\,I~Nicolaescu},  \emph{Notes on Seiberg--Witten theory}, American
Mathematical Society, Providence RI (2000)

\bibitem{pao:torus-I} {\bf P\,S Pao}, \emph{The topological structure
of $4$--manifolds with effective torus actions. {I}},
Trans. Amer. Math. Soc. {227} (1977) 279--317

\bibitem{quinn:bordism} {\bf F Quinn}, \emph{Open book decompositions,
and the bordism of automorphisms}, Topology, {18} (1979) 55--73


\bibitem{rosenberg-stolz:psc} {\bf J Rosenberg}, {\bf S Stolz},
\emph{Metrics of positive scalar curvature and connections with
surgery}, from: ``Surveys on surgery theory, Vol. 2'', Princeton
Univ. Press, Princeton, NJ (2001) 353--386

\bibitem{ruberman:isotopy} {\bf D Ruberman}, \emph{An obstruction to
smooth isotopy in dimension $4$}, Math. Res. Lett. {5} (1998) 743--758,
{\tt arxiv:math.GT/9807041}

\bibitem{ruberman:polyisotopy} {\bf D Ruberman}, \emph{A polynomial
invariant of diffeomorphisms of 4-manifolds}, from: ``Proceedings of
the Kirbyfest (Berkeley, CA, 1998)'', Geom. Topol. Monogr. 2 (1999)
473--488, {\tt arxiv:math.GT/9911260}
 
\bibitem{schoen-yau:psc} {\bf R~Schoen}, {\bf S\,T Yau}, \emph{On the
structure of manifolds with positive scalar curvature}, Manuscripta
Math. {28} (1979) 159--183

\bibitem{taubes:l2} {\bf C~Taubes}, \emph{{$L^2$}--moduli spaces on
$4$---manifolds with cylindrical ends}, Monographs in Geometry and
Topology {1}, International Press (1994)

\end{thebibliography}
\end{document}